\documentclass[11pt]{amsart}
\usepackage{lmodern}
\usepackage{amssymb,xcolor,array}
\usepackage[numbers,square,sort&compress]{natbib}
\usepackage[hidelinks]{hyperref}
\hypersetup{
  pdftitle={Central fibers of quantum skew-symmetric matrices at roots of unity},
  pdfauthor={Xiaomin Tang and Yu Zhang}}
\numberwithin{equation}{section}
\theoremstyle{plain}
\newtheorem{theorem}{Theorem}[section]
\newtheorem{corollary}[theorem]{Corollary}
\newtheorem{proposition}[theorem]{Proposition}

\newtheorem*{boundarytheorem}{Boundary theorem}
\theoremstyle{definition}
\newtheorem{example}[theorem]{Example}
\theoremstyle{remark}
\newtheorem{remark}[theorem]{Remark}
\allowdisplaybreaks[2]
\begin{document}
\title[Central fibers of quantum skew-symmetric matrices]
 {Central fibers of quantum skew-symmetric matrices at roots of unity}
\author[X. Tang]{Xiaomin Tang}
\address{School of Mathematics and Statistics, Heilongjiang University,
 Harbin 150080, China}
\email{tangxm@hlju.edu.cn}
\author[Y. Zhang]{Yu Zhang}
\address{School of Mathematics and Statistics, Heilongjiang University,
 Harbin 150080, China}
\email{zhangyu@hlju.edu.cn}
\thanks{Corresponding author: Xiaomin Tang
 (ORCID \href{https://orcid.org/0000-0002-9054-1286}{0000-0002-9054-1286}).}

\begin{abstract}
For the quantum coordinate algebra of skew-symmetric matrices at a root
of unity of odd order, we study the Azumaya locus over the full center
and its relation to quantum cluster seed charts. In every even size at
least six, we construct an Azumaya open set on which the points missed
by every seed chart form a torus. The Brauer index at the generic point
of this torus is the square root of the generic index on the open set,
although the closed fibers on the entire open set are matrix algebras
of the same size. The calculation uses a simultaneous exchange
factorization into rank-one algebras and a residual quantum torus.
We also describe the full center over the power center and compute
simple modules and projective radical layers in seed charts and in
specified boundary fibers. An independent calculation compares the
order-three endpoint algebra with two subregular Harish--Chandra
reductions in type \(A_2\).
\end{abstract}
\maketitle

\section{Introduction}\label{sec:intro}

At a root of unity, a quantum cluster chart gives useful coordinates
for studying representations of an algebra over its center. The chart
inverts the central powers of the mutable variables, so it cannot see
a central character at which some mutable variable in every seed
fails to be invertible. This raises a question: can a central character
missed by all seed charts still be Azumaya? For quantum skew-symmetric
matrices we give an explicit family of such points and compute the
Brauer index along it.

Let \(A=A_\zeta(n)\) be the quantum coordinate algebra of
\[
 \operatorname{Skew}_n
  =\{X\in\operatorname{Mat}_n(\mathbb C):X^{\mathsf T}=-X\},
 \qquad n\ge4,
\]
where \(\zeta\) has odd order \(\ell>1\).
The defining relations are given in \eqref{eq:skew-relations}.
At generic \(q\), we use the explicit quantum seed of \citet{ZT}
for this algebra. For the invariant-theoretic and reflection-equation
background, see \citep{Strickland,JR,JZ}.
The specialized algebra is finite over its center \(Z(A)\).
A character \(\chi:Z(A)\to\mathbb C\) determines the full-center
fiber
\[
 A_\chi=A\otimes_{Z(A),\chi}\mathbb C .
\]
The Azumaya locus consists of the points where \(A\) is locally a
matrix algebra in the \'etale topology of \(\operatorname{Spec}Z(A)\).

A quantum seed has quasi-commuting variables \(Y_i\), an integral
commutation matrix \(\Lambda\), and an exchange matrix. Mutation
replaces one mutable variable by the sum of two normalized Laurent
monomials. Frozen variables are not replaced, and we do not impose
their invertibility. At \(\zeta\), the powers \(Y_i^\ell\) are central
\citep{NTY,HLY}. The seed chart in the full center is
\[
 U_\Sigma=\bigcap_{i\in I_{\mathrm{mut}}(\Sigma)}D(Y_i^\ell)
 \subseteq\operatorname{Spec}Z(A),
 \qquad D(g)=\{\mathfrak p:g\notin\mathfrak p\}.
\]
For \(\chi\in U_\Sigma\), the mutable variables have inverses
\(Y_i^{\ell-1}/\chi(Y_i^\ell)\) in \(A_\chi\). We will identify an
open set on which failure of this invertibility test in every seed
does not obstruct the Azumaya property.

Over a field, a central simple algebra has the form \(M_s(D')\) for
a central division algebra \(D'\). Its degree is \(s\deg D'\), and
its Brauer index is \(\deg D'\). Thus identical closed matrix fibers
do not imply that the index is constant at generic points of
subvarieties. The following theorem makes the difference explicit.

Suppose \(n=2m\ge6\), and set
\[
 \mathcal I=\{(i,j):1\le i<j\le 2m-1,\ i,j\text{ odd}\},
 \qquad r=|\mathcal I|=\binom m2 .
\]
For \(a\in\mathcal I\), let \(X_a\) be the initial variable at \(a\)
and \(W_a\) its mutation. Let \(U\subseteq\operatorname{Spec}Z(A)\)
be the principal open where all initial powers indexed outside
\(\mathcal I\) are units, and set
\[
 \Delta=\left(U\setminus\bigcup_\Sigma U_\Sigma\right)_{\mathrm{red}}.
\]
The main boundary result combines Theorem~\ref{thm:odd-pair-open}
and Corollary~\ref{cor:odd-pair-brauer}.

\begin{boundarytheorem}
The restriction of \(A\) to \(U\) is Azumaya. On this open one has
\[
 \Delta=V_U\!\left(X_a^\ell,W_a^\ell:a\in\mathcal I\right)_{\mathrm{red}}
       \cong (\mathbb C^\times)^{m^2}.
\]
If \(\eta_U\) and \(\eta_\Delta\) are the generic points of \(U\)
and \(\Delta\), respectively, then
\[
 \operatorname{ind}(A_{\eta_U})=\ell^{2r},\qquad
 \operatorname{ind}(A_{\eta_\Delta})=\ell^r,\qquad
 A_\chi\cong M_{\ell^{2r}}(\mathbb C)
 \quad(\chi\in U(\mathbb C)).
\]
\end{boundarytheorem}

The proof uses simultaneous exchanges at the positions in
\(\mathcal I\). On \(\Delta\), each rank-one factor becomes a split
matrix algebra; the residual quantum torus supplies the division
algebra at \(\eta_\Delta\). The two indices above are taken over
function fields of loci in the \emph{same full center}.
The change of center is a separate statement: the restriction of
\(\operatorname{Spec} Z(A)\to\operatorname{Spec} Z_0\), where
\(Z_0=\mathbb C[x_{ij}^{\ell}:i<j]\), maps \(\Delta\) onto its
reduced image by a finite \'etale map of degree \(\ell^m\).

\citet{BG02} study central fibers of quantized
function algebras of semisimple groups. \citet{BY}
relate Azumaya loci of PI algebras to discriminant ideals. In the
cluster setting, the power subalgebra of Nguyen, Trampel and
Yakimov \cite{NTY} provides the commutative coordinates, and
\citet[Theorem~6.1 and Section~6.3]{MNTY} bound the fully Azumaya
locus over that subalgebra under their hypotheses. Our calculations
use the full center. They determine the radical layers at
simultaneous frozen zeros, and, on an even-size open, identify the
points missed by every seed and compute their generic Brauer index.
Bera and Mukherjee \citep{BeraMukherjee} determine the center and
classify simple modules for a two-parameter quantum matrix algebra
at roots of unity. Our algebra is the skew-matrix algebra; the
boundary calculation above uses its exchange relations to determine
a Brauer class where no seed chart applies.

For a full-center character \(\chi\) and its restriction
\(\chi_0=\chi|_{Z_0}\), the natural surjection
\[
 A\otimes_{Z_0,\chi_0}\mathbb C\ \,\longrightarrow\,\,
 A\otimes_{Z(A),\chi}\mathbb C
\]
imposes the additional central relations. We determine them first,
then use them to compute projective layers on seed charts and at
points where the two exchange-coordinate powers vanish.

For the fibers that are not Azumaya, we also ask what can be
recovered about projective covers from exchange coordinates.
Let \(P_\alpha\) be the projective cover of a simple \(A_\chi\)-module
\(L_\alpha\), and let \(J=\operatorname{rad}(A_\chi)\).
The two
multiplicity invariants of interest are
\begin{equation}\label{eq:intro-multiplicities}
 C_{\beta\alpha}(t)
   =\sum_{d\ge0}[J^dP_\alpha/J^{d+1}P_\alpha:L_\beta]t^d,
 \qquad
 [P_\alpha:L_\beta]=C_{\beta\alpha}(1).
\end{equation}
Here \([V:L]\) denotes composition multiplicity. We ask how much
of \eqref{eq:intro-multiplicities} can be computed from cluster
coordinates when a central character sets some coordinates to zero.

This question is motivated by Tanisaki's theorem
\citep[Theorem~10.4]{Tanisaki} on Lusztig's multiplicity conjecture
for the De Concini--Kac quantum group at roots of unity. In
Tanisaki's setting one fixes a Frobenius character and works
with a generalized Harish--Chandra character. His localization
identifies the simples \(L_\tau\) and projective covers
\(E_\sigma\) with objects in an exotic heart on the geometric side.
The associated equivariant \(K\)-theory expansion gives
\begin{equation}\label{eq:tanisaki-motivation}
 [\mathcal E'_\sigma]=\sum_\tau n_{\tau,\sigma}(v)
       [\mathcal L'_\tau]
 \quad\Longrightarrow\quad
 [E_\sigma:L_\tau]=n_{\tau,\sigma}(1).
\end{equation}
We ask the corresponding multiplicity question for quantum
skew-symmetric matrices, using cluster charts and exchange relations.
The variable \(t\) in \eqref{eq:intro-multiplicities} records radical
degree, not the equivariant grading in \eqref{eq:tanisaki-motivation}.

Write \(N=\binom n2\), and index the specialized initial seed by
\(Y_{r,c}\) with \(r<c\). We call the set of frozen coordinates
whose central powers vanish the frozen zero pattern.
Put \(F_r=Y_{r,n}\) for \(1\le r\le n-2\).
In even size \(n=2m\), write \(H,B\) for the two terminal
coordinates of the adjacency chains on \([1,n]\) and
\([2,n-1]\), respectively. In odd size \(n=2m+1\), write \(B_L,B_R\)
for the corresponding terminal coordinates on the intervals
\([1,n-1]\), \([2,n]\). The construction of these coordinates
and the normalization of their monomials are given in
Section~\ref{sec:integral-structure}.
The next result identifies the additional central parameters:
they are indexed by a radical supported entirely on the frozen
coordinates.

\begin{theorem}\label{thm:intro-center}
Let \(K=\ker(\Lambda\bmod\ell)\). In the initial seed the
following are bases over \(\mathbb Z/\ell\mathbb Z\):
\begin{align*}
 n=2m:quad&
 e_H,\quad e_B-e_{F_1},\quad
 e_{F_{n-r}}-e_{F_r}\quad(2\le r\le m-1);\\
 n=2m+1:\quad&
 e_{B_L}+e_{B_R}-e_{F_1},\quad
 e_{F_{n-r}}-e_{F_r}\quad(2\le r\le m).
\end{align*}
For \(k\in K\), let \(\bar k\) be its representative with all
entries in \([0,\ell-1]\), and let \(Y^{\bar k}\) be the based
frozen monomial with that exponent. Then
\begin{equation}\label{eq:main-center}
 Z(A)=\bigoplus_{k\in K} Z_0 Y^{\bar k},
 \qquad
 \operatorname{rank}_{Z_0}Z(A)=\ell^{\lfloor n/2\rfloor}.
\end{equation}
\end{theorem}

Let \(S\) be the frozen zero set of a full-center character
in a seed chart, \(U\) the set of all remaining coordinates
and \(\pi_S\) the coordinate projection to \(S\).
Write \(K_U=\ker(\Lambda_{UU}\bmod\ell)\),
\(K_0=\ker(\pi_S|_K)\), and put
\[
 \mathcal D_S=
 \{d\in\mathbb N^S:\ d\not\ge e
   \text{ for every }0\ne e\in\mathbb N^S
   \text{ with }e\bmod\ell\in\pi_S(K)\},
\]
with componentwise order. We call the full-center character
Azumaya if the algebra over the center is Azumaya at that point.
In \(\mathcal D_S\) we retain the zero-coordinate monomials not
annihilated by the central relations. Their total degrees measure
how far the corresponding factors survive in the radical filtration.

\begin{theorem}\label{thm:intro-chart}
For a full-center character in a seed chart, the simple modules
and each projective cover satisfy
\[
 \#\operatorname{Irr}(A_\chi)=\frac{|K_U|}{|K_0|},\qquad
 \dim L=\sqrt{\frac{\ell^{|U|}}{|K_U|}},\qquad
 \sum_\beta C_{\beta\alpha}(t)=\sum_{d\in\mathcal D_S}t^{|d|}.
\]
The last expression is a product over the supports in
Theorem~\ref{thm:intro-center} and, in even size, the unpaired
frozen coordinate \(F_m\). The \(\chi\) on the Azumaya locus
are exactly those with at most one zero on each support and
with \(F_m\ne0\) in even size. The numbers of these frozen zero
patterns are \(2\cdot3^{m-1}\) for \(n=2m\) and
\(4\cdot3^{m-1}\) for \(n=2m+1\).
\end{theorem}

Corollary~\ref{cor:global-azumaya} passes from the chart criterion
to the whole center: in odd size the Azumaya and regular loci
coincide, and in even size the Azumaya locus is the regular locus
with the divisor \(F_m^\ell=0\) removed.

The boundary results distinguish two kinds of cancellation. Write
\(XW=A_++A_-\) for an exchange and call \(X^\ell=W^\ell=0\) a
collision. When both exchange monomials are units, the rank-one
factorization in Theorem~\ref{thm:unit-exchange-patch} identifies
the collision fiber with a fiber on a seed chart, possibly at a
different point. When all other initial mutable powers are units
and \(n\ge6\), Proposition~\ref{prop:single-mutable-coverage} sharpens
this to a criterion for the point itself:
\[
 \begin{array}{c|l}
 \text{exchange monomials at a collision} & \text{coverage}\\ \hline
 A_+^\ell A_-^\ell\ne0 & \text{a chart after two mutations}\\
 A_+^\ell=A_-^\ell=0 & \text{no seed chart.}
 \end{array}
\]
The second line gives exactly the three frozen configurations in
that proposition. Theorem~\ref{thm:all-single-nodal-layers} computes
their labelled layers. In the first-row terminal case with only
\(F_1,F_2\) zero, the polynomial is \(1+2\sum_{d=1}^{\ell-1}t^d\),
rather than \(1\) on a chart with the same frozen zeros.

Section~\ref{sec:block} gives a contrasting block-diagonal family
and computes its labelled radical layers over \(Z_0\) by a
cut filtration on the PBW basis.

The comparison with Tanisaki uses the same doubled cycle on
both sides: for the order-three skew endpoint it comes from
the exchange pair, and for the subregular \(A_2\) block from
contracting its invertible bonds. Proposition~\ref{prop:endpoint-thickening}
identifies the generalized and exact reductions with two quotients
of the completed endpoint, up to Morita equivalence. Their Cartan
matrices are \(6\mathbf 1_{3\times3}\) and
\(2\mathbf 1_{3\times3}-I_3\), respectively. The comparison is
of these central reductions; it does not identify Tanisaki's
geometric categories or equivariant grading.

\section{Preliminaries}\label{sec:preliminaries}

\subsection{The quantum skew-matrix algebra}\label{subsec:skew-definition}
Let \(R=\mathbb C[v^{\pm1}]\), \(q=v^2\), and
\(\mathcal I_n=\{(i,j):1\le i<j\le n\}\), with lexicographic
order. A classical skew-symmetric matrix satisfies
\(X^{\mathsf T}=-X\); its entries above the diagonal give
\(\mathbb C[\operatorname{Skew}_n]=\mathbb C[x_{ij}:(i,j)\in\mathcal I_n]\).
The quantum coordinate algebra used here is the associative algebra
\begin{equation}\label{eq:skew-presentation}
 A_R(n)=R\langle x_{ij}:(i,j)\in\mathcal I_n\rangle/
             (\text{relations \eqref{eq:skew-relations}}).
\end{equation}
Put \(\eta=q-q^{-1}\). In the first line of
\eqref{eq:skew-relations}, \(a<b\) are pairs sharing an index;
in the other three lines \(r<s<t<u\):
\begin{equation}\label{eq:skew-relations}
\begin{aligned}
 x_bx_a       &=q^{-1}x_ax_b,\\
 x_{st}x_{ru} &=x_{ru}x_{st},\\
 x_{su}x_{rt} &=x_{rt}x_{su}-\eta x_{ru}x_{st},\\
 x_{tu}x_{rs} &=x_{rs}x_{tu}-\eta(x_{rt}x_{su}-q x_{ru}x_{st}).
\end{aligned}
\end{equation}
The four lines cover every pair of distinct generators: a shared
index, then nested, crossing and separated pairs. Thus
\(A_R(n)/(v-1)A_R(n)=\mathbb C[\operatorname{Skew}_n]\).
We use the upper-entry convention of \citep{Strickland}; the
noncommutative algebra is specified by \eqref{eq:skew-relations},
not by imposing a transpose identity on a matrix of generators.

The generic and specialized algebras will be
\begin{equation}\label{eq:generic-specialized}
 A_{\mathbb C(v)}=A_R(n)\otimes_R\mathbb C(v),\qquad
 A_\zeta(n)=A_R(n)/(v-\nu)A_R(n),\qquad
 \nu=\zeta^{(\ell+1)/2}.
\end{equation}
Here \(\zeta\) is primitive of odd order \(\ell>1\), so that
\(\nu^2=\zeta\) and \(\nu^\ell=1\). We abbreviate \(A_R(n)\)
to \(A_R\) when \(n\) is fixed.

For an increasing even set \(I=(i_1,\ldots,i_{2h})\), the
quantum Pfaffian in this convention is
\begin{equation}\label{eq:qpfaffian}
 [\varnothing]_q=1,\qquad
 [I]_q=\sum_{j=2}^{2h}(-q)^{j-2}x_{i_1i_j}
              [I\setminus\{i_1,i_j\}]_q.
\end{equation}
The factor order is part of the definition. In particular,
\[
 [r,s,t,u]_q
   =x_{rs}x_{tu}-q x_{rt}x_{su}+q^2x_{ru}x_{st}.
\]
This is the first-row convention of \citep[(1.6)]{Strickland}
and specializes to the ordinary Pfaffian at \(q=1\).

\subsection{Quantum tori and seed mutation}\label{subsec:seed-prelim}
A based quantum torus \(\mathcal T_\Lambda\) for an integral
alternating form \(\Lambda\) has basis \(M(a)\), \(a\in \mathbb Z^N\).
The multiplication and its expression in ordered generators are
\[
  \begin{aligned}
 M(a)M(b)&=v^{a^{\mathsf T}\Lambda b}M(a+b),\\
 M(a)&=v^{-\sum_{i<j}\Lambda_{ij}a_ia_j}
                 Y_1^{a_1}\cdots Y_N^{a_N},\qquad Y_i=M(e_i).
 \end{aligned}
\]
In particular,
\(Y_iY_j=q^{\Lambda_{ij}}Y_jY_i\). A quantum seed consists of
a based torus embedded in the ambient skew field, a set \(I_{\rm mut}\)
of mutable indices, and an integer matrix \(\widetilde B\) whose
columns \(b_k\) are indexed by \(I_{\rm mut}\). The principal part
on \(I_{\rm mut}\) is skew-symmetrizable. In the compatible sign
convention the seed data satisfy
\[
 b_k^{\mathsf T}\Lambda=d_k e_k^{\mathsf T},\qquad
 d_k\in\mathbb Z_{>0}\quad(k\in I_{\rm mut}).
\]
The other indices \(I_{\rm fr}\) are frozen and cannot be mutated.
For \([b]_+\) the coordinatewise positive part, mutation at \(k\)
fixes the other coordinates and replaces \(Y_k\) by
\[
  Y'_k=M(-e_k+[b_k]_+)+M(-e_k+[-b_k]_+).
\]
The new exchange matrix \(\widetilde B'=(b'_{ij})\) is
\[
 b'_{ij}=
 \begin{cases}
 -b_{ij},&i=k\ \text{or}\ j=k,\\
 b_{ij}+[b_{ik}]_+[b_{kj}]_+-[-b_{ik}]_+[-b_{kj}]_+,&i,j\ne k.
 \end{cases}
\]
The new commutation form is \(E_k^{\mathsf T}\Lambda E_k\),
where \(E_k\) fixes the other coordinate vectors and
\(E_ke_k=-e_k+[b_k]_+\). See \citep{BZ}.
The quantum cluster algebra is the subalgebra generated by
the seed coordinates of all seeds obtained by successive mutations.
We do not adjoin frozen inverses. For \(S\subseteq I_{\rm fr}\)
the \emph{mixed torus} \(T_S\) is the span of the \(M(a)\) with
\(a_i\ge0\) for \(i\in S\). The \emph{quantum Laurent phenomenon}
is the inclusion in these Laurent rings with only the permitted
denominators; in particular no frozen exponent is negative in
the mixed ring \(T_{I_{\rm fr}}\) of a seed. For the specialized
algebra, a seed chart is the central open set
\(\bigcap_{i\in I_{\rm mut}}D(Y_i^\ell)\). Section~\ref{sec:integral-structure}
establishes the localization identity \eqref{eq:all-mixed-charts}
and not merely a Laurent expansion in its fraction field.

\subsection{From skew polynomial extensions to quantum seeds}
\label{subsec:cgl-prelim}

The seed used below is produced before the specialization. We
recall the relevant construction and separate it from its
integral and root-of-unity steps.

Let \(K=\mathbb C(v)\), \(q=v^2\). An Ore extension
\(B[x;\theta,\delta]\) is the polynomial extension as a left
\(B\)-module with multiplication
\[
 xb=\theta(b)x+\delta(b),\qquad
 \delta(bb')=\theta(b)\delta(b')+\delta(b)b'.
\]
Here \(\theta\) is an automorphism and \(\delta\) is a
\(\theta\)-skew derivation. An iterated Ore extension
\[
 K[x_1][x_2;\theta_2,\delta_2]\cdots[x_N;\theta_N,\delta_N]
\]
is called a CGL extension if an algebraic torus acts rationally
by automorphisms, each \(x_i\) is a torus eigenvector, each
\(\delta_i\) is locally nilpotent on the preceding algebra,
and \(\theta_i\) is induced by a torus element whose eigenvalue
on \(x_i\) is not a root of unity. The symmetric condition asks
that \(\delta_j(x_i)\) belong to the subalgebra on
\(x_{i+1},\ldots,x_{j-1}\) for \(i<j\) and that the inverse
commutation scalars likewise be induced by torus elements in
the reverse order. Symmetry is about the two orders of the
extension, not about symmetry of the quantum matrix.

A homogeneous prime element is a torus eigenvector \(y\)
with \(yB=By\) and \(B/(y)\) a domain. The CGL prime theorem
\citep[Theorem~3.6]{GY} gives monic prefix primes \(y_j\):
\[
 y_j=
 \begin{cases}
 x_j,&\text{if no predecessor exists},\\
 y_{p(j)}x_j-c_j,&p(j)<j,\quad c_j\text{ in the preceding prefix}.
 \end{cases}
\]
The prefix primes, up to scalars, are the last chain members
present in that prefix. The leading monomial of \(y_j\)
records its predecessor chain.

Under the square-root, reverse-eigenvalue and interval
normalization hypotheses of \citet[Theorem~8.2 and
Proposition~8.9]{GY}, the normalized \(y_j\) form a seed
with last chain members frozen. The interval seeds in its
mutation class contain every \(x_i\) up to normalization.
Section~\ref{sec:integral-structure} checks these hypotheses
and proves that the seed localizations survive specialization.

\subsection{Central powers, fibers and radical layers}\label{subsec:fiber-prelim}
For a based quantum torus at \(v=\nu\), write
\(K_\Lambda=\ker(\Lambda\bmod\ell)\) in \((\mathbb Z/\ell)^N\).
Conjugation by the \(Y_i\) is diagonal on the monomial basis,
so
\begin{equation}\label{eq:torus-central-test}
 M(a)\in Z(\mathcal T_\Lambda)
 \quad\Longleftrightarrow\quad
 \Lambda a\equiv0\bmod\ell,\qquad
 M(a)^\ell=M(\ell a)=\prod_i(Y_i^\ell)^{a_i}.
\end{equation}
For the mixed torus one retains the same central exponents
subject to the nonnegativity restrictions. In particular the
coordinate powers are central. If a compatible exchange has
diagonal entry \(d_k\) prime to \(\ell\), its two monomials
\(T_+,T_-\) satisfy
\(T_+T_-=\zeta^{\pm d_k}T_-T_+\). The quantum binomial
identity then specializes to
\[
 (T_++T_-)^\ell=T_+^\ell+T_-^\ell.
\]
Thus the \(\ell\)-th powers satisfy the commutative exchange.
For the root-of-unity seed algebras considered here the
central embedding of \citep[Theorem~4.6]{NTY} identifies the
commutative cluster algebra with the subalgebra generated by
these powers, called the \emph{Frobenius subalgebra}.
In Section~\ref{sec:integral-structure} we verify the specialization
hypotheses for \(A_\zeta(n)\); the resulting subalgebra is
\(Z_0=\mathbb C[x_{ij}^\ell:i<j]\).

For \(Z_0\subseteq Z(D)\) with \(D\) finite over \(Z_0\), the
power and full-center fibers at compatible characters
\(\chi_0\) and \(\chi\) are
\begin{equation}\label{eq:two-fibers}
 D_{\chi_0}=D\otimes_{Z_0,\chi_0}\mathbb C,\qquad
 D_\chi=D\otimes_{Z(D),\chi}\mathbb C,\qquad
 D_{\chi_0}\twoheadrightarrow D_\chi
 \quad(\chi|_{Z_0}=\chi_0).
\end{equation}
Only the second fiber imposes the values of all central elements.
We call a point Azumaya if the algebra at the corresponding
local ring of the full center is Azumaya: for an algebra \(B\)
over a commutative ring \(C\) this means \(B\) is a faithful
finite projective \(C\)-module and
\[
 B\otimes_C B^{\rm op}\longrightarrow\operatorname{End}_C(B),
 \qquad a\otimes b\longmapsto(d\mapsto adb)
\]
is an isomorphism. A complex fiber at such a point is a
matrix algebra.

For a finite-dimensional algebra \(B\), its Jacobson radical is
\(J=\operatorname{Jac}(B)=\bigcap_{L\ {\rm simple}}\operatorname{Ann}_B(L)\).
It is nilpotent and \(B/J\) is semisimple. Choose a primitive
idempotent \(e_\alpha\) for each simple label, so that
\(P_\alpha=Be_\alpha\) and \(L_\alpha=P_\alpha/JP_\alpha\).
For \(e=\sum_\alpha e_\alpha\) one has \(BeB=B\), and \(eBe\) is
the basic algebra, with one-dimensional simple modules.
The equivalence \(M\mapsto eM\) preserves the radical-layer
polynomials \eqref{eq:intro-multiplicities}.

\subsection{A two-generator exchange and its signed degrees}\label{subsec:signed-prelim}
The exchange computations use the normal forms associated
with rank-one generalized Weyl algebras \citep{Bavula}.
For an automorphism \(\sigma\) of an algebra \(D\), the skew Laurent
ring \(D[X^{\pm1};\sigma]\) consists of finite sums
\(\sum_i d_iX^i\) with \(Xd=\sigma(d)X\).
If \(XW=h\in D\) in an overring with \(X\) invertible, then
\[
 W=X^{-1}h,\qquad WX=\sigma^{-1}(h),\qquad
 Wd=\sigma^{-1}(hd)\,X^{-1}.
\]
In the applications below \(h\) is a sum of two monomials
normalizing \(D\) in the same way, or lies in a commutative
coefficient algebra. Thus both \(X\) and \(W\) can be moved to
the right of coefficients. The grading by
\(\deg D=0,\ \deg X=1,\ \deg W=-1\) is called the
\emph{signed degree}. Removing \(XW\) or \(WX\) reduces a
homogeneous word to a coefficient times a power of \(X\)
or a power of \(W\). Independence requires a separate argument.
The signed degree is distinct from the radical degree.

\section{Integral seeds and their specializations}\label{sec:integral-structure}

We now verify the generic seed, integral form and specialized
chart identities for the presentation \eqref{eq:skew-presentation}.

\subsection{The integral presentation and the CGL seed}

\noindent\textbf{Ordered PBW basis.}\enspace
Write the upper entries as \(x_1,\ldots,x_N\) in lexicographic order.
An ordered PBW (Poincar\'e--Birkhoff--Witt) basis here means
\(\{x_1^{a_1}\cdots x_N^{a_N}:a\in\mathbb N^N\}\).
For a nonzero polynomial in this basis, the leading exponent
is the maximal exponent for comparison starting at \(a_N\).
Every nonscalar correction in the straightening rules decreases
the comparison of the exponent multisets starting at the largest index;
a scalar interchange reduces the inversion count if that multiset
is unchanged. Thus the ordered monomials
\(x_1^{a_1}\cdots x_N^{a_N}\) span over \(R\).
Strickland's generic ordered PBW theorem
\citep[Proposition~1.1]{Strickland}, in the same upper-entry
convention (see also her proof of Theorem~1.5 for the integral
form), gives their independence over \(\mathbb C(v)\). A relation
over \(R\) would also be a relation over that field. The ordered
monomials are therefore an \(R\)-basis and remain a basis under
specialization.

For \(a=(r,c)\), denote the earlier subalgebra by \(A_{<a}\).
Let \(\nu(a)=e_r+e_c\), and use the Ore convention
\(x_a b=\theta_a(b)x_a+\delta_a(b)\) for \(b\in A_{<a}\).
The ordered relations give
\[
 \deg_{\rm row}x_a=\nu(a),\qquad
 \theta_a(x_b)=q^{-(\nu(a),\nu(b))}x_b,\qquad
 \theta_a(x_a)=q^{-2}x_a
\]
where the last equality uses the torus action extending \(\theta_a\).
Moreover,
\[
 \delta_a(x_b)\in \mathbb C(v)\langle x_j:b<j<a\rangle,\qquad
 \delta_a^2(x_b)=0 .
\]
In the second equality every correction factor shares an endpoint
with \(a\), hence is killed by \(\delta_a\). The skew Leibniz rule,
with \(\theta_a\) diagonal, gives
\[
 \delta_a^{h+1}(x_{b_1}\cdots x_{b_h})=0 .
\]
Reversing the order negates the torus exponents and preserves
the correction bounds. This gives a symmetric CGL extension,
with reverse eigenvalues \(q^2\).

The predecessor map is
\begin{equation}\label{eq:pred-map}
 p(r,c)=
 \begin{cases}
  (r-1,c-1),&r>1,\ c-r\ge2,\\
  (r-2,r-1),&c=r+1,\ r\ge3,\\
  -\infty,&r=1\ \text{or}\ (r,c)=(2,3).
 \end{cases}
\end{equation}
The two Ore corrections at the proposed predecessor,
read directly from \eqref{eq:skew-relations}, are
\[
\begin{aligned}
 \delta_{(r,c)}(x_{r-1,c-1})
   &=-\eta\,x_{r-1,c}x_{r,c-1} &&(c-r\ge2),\\
 \delta_{(r,r+1)}(x_{r-2,r-1})
   &=-\eta\bigl(x_{r-2,r}x_{r-1,r+1}
                  -q x_{r-2,r+1}x_{r-1,r}\bigr).
\end{aligned}
\]
Here is the lead check in the induction on \(a=(r,c)\).
If \(p=p(a)\ne-\infty\), it is terminal in the earlier prefix
by the induction hypothesis. Write
\(Q_p=\operatorname{lt}(y_p)/x_p\) for the
ordered product of the other chain entries, as a commutative
PBW lead monomial. In the two cases of \eqref{eq:pred-map}
the Ore correction of \(y_p\) has a nonzero term whose PBW monomial is listed below.
Table~\ref{tab:lead} records the ordered monomials.
\begin{table}[htbp]
\caption{Leading terms detecting the two predecessors.}\label{tab:lead}
\centering
\begin{tabular}{>{$}c<{$}|>{$}c<{$}|>{$}c<{$}}
 p & \text{last vertex in the row weight of }y_p
     & \text{lead witness in }\delta_a(y_p)\\ \hline
 (r-1,c-1)&c-1&Q_p\,x_{r-1,c}x_{r,c-1}\\
 (r-2,r-1)&r-1&Q_p\,x_{r-2,r}x_{r-1,r+1}
 \end{tabular}
\end{table}
The displayed products denote exponent multisets.
For the first row, vertex \(c-1\) occurs once and vertex \(c\)
does not occur. A term containing \(x_{r,c-1}\) must originate
at the unique factor \((b,c-1)\). If \(b<r-1\), the new
\(x_{b,c}\) is smaller than \(x_{r-1,c}\) in the PBW
comparison; if \(b=r-1\), the substitution preserves the order
on residual exponent multisets. No straightening term can
create a larger factor in row \(r\). For the second row both
\(r,r+1\) are absent. Among corrections at an entry \((b,d)\)
in \(y_p\), maximize first the new factor \(x_{d,r+1}\)
and then \(x_{b,r}\). The maximum requires
\((b,d)=(r-2,r-1)=p\); on the residual multisets it again
preserves order. In each row the coefficient is a nonzero
Laurent unit times \(q-q^{-1}\). Thus \(\delta_a(y_p)\ne0\).
If \(y_p\) had remained normal in the new prefix, homogeneity
and the one-dimensional weight space of \(x_a\) would give
\(x_a y_p=\lambda y_p x_a\); lead comparison forces the
\(\theta_a\)-scalar, in contradiction to \(\delta_a(y_p)\ne0\).
The CGL prime theorem therefore identifies \(p\) as the
predecessor. When \(p=-\infty\), all Ore corrections vanish.

We next determine the interval primes needed for the seed
normalization. For \(r<s<t<u\), consider either of the endpoint
pairs \((L,U)=(x_{rt},x_{su})\) with \(s=r+1,\ u=t+1\),
or \((L,U)=(x_{rs},x_{tu})\) with
\((s,t,u)=(r+1,r+2,r+3)\).
Put \(N_0=x_{ru}x_{st}\), and in the second case put
\(C_0=x_{rt}x_{su}\).
The interval consists of the lexicographic generators from \(L\)
to \(U\), inclusive. Before \(U\), each generator either shares
an index with \(L\) or is nested with it. Thus \(L\) is normal
in the prefix without \(U\), whereas
\[
 [L,U]=\begin{cases}
 \eta N_0,&(L,U)=(x_{rt},x_{su}),\\
 \eta(C_0-qN_0),&(L,U)=(x_{rs},x_{tu})
 \end{cases}
\]
is nonzero. By the CGL prime theorem the monic prime at \(U\)
in this interval has predecessor \(L\).
Its row weight uses \(r,s,t,u\) once each. The possible
monomials of that weight give the following two forms:
\[
 Q=LU+\beta N_0,\qquad Q=LU+\alpha C_0+\beta N_0,
\]
respectively. In the first case the remaining pairing uses
\(x_{rs}\), which occurs before the start of the interval.
Normality and leading terms give \([L,Q]=0\), since \(L,U\)
have disjoint endpoints. Also \(LN_0=q^2 N_0 L\) and,
where applicable, \(LC_0=q^2 C_0 L\). Consequently
\[
\begin{aligned}
 \relax[L,LU+\beta N_0]
 &=L\bigl(\eta+(1-q^{-2})\beta\bigr)N_0,\\
 [L,LU+\alpha C_0+\beta N_0]
 &=L\bigl((\eta+(1-q^{-2})\alpha)C_0
        +(-q\eta+(1-q^{-2})\beta)N_0\bigr).
\end{aligned}
\]
The PBW basis separates \(C_0,N_0\), so
\(\eta/(1-q^{-2})=q\) yields
\[
 Q=LU-qN_0,\qquad Q=LU-qC_0+q^2N_0,
\]
respectively. Thus the coefficient identities follow from
normality of the existing CGL prime, without a finite check.

In either case the leading term of \(LU-Q\) is \(q x_i x_j\),
where \(i<j\) have disjoint endpoints and each shares one
endpoint with \(L=x_a\), \(a<i<j\). In terms of the
straightening pairing \(M\) (with \(M_{ij}=1\) for \(i<j\)
sharing an index and zero for disjoint pairs), this gives
\[
 M_{ai}=M_{aj}=1,\qquad M_{ij}=0,\qquad
 -\tfrac12(-M_{ai}-M_{aj}+M_{ij})=1 .
\]
Hence the required coefficient is
\(S_\nu(-e_a+e_i+e_j)=q\), using the GY square-root
normalization \(S_\nu\) of \citep[condition (6.13)]{GY}.
The square-root scalars have no generic \(2\)-torsion.
Thus the hypotheses of \citep[Theorem~8.2 and
Proposition~8.9]{GY} apply: the generic algebra is the
quantum cluster algebra without frozen inverses, and each
\(x_i\) occurs up to a Laurent unit in an interval seed.
The frozen variables of the initial seed are the chain
terminals.

For instance, in size four the monic coordinates in lexicographic
order are
\[
 x_{12},\quad x_{13},\quad x_{14},\quad x_{23},\quad
 x_{13}x_{24}-q x_{14}x_{23},\quad [1,2,3,4]_q.
\]
The first two are mutable and the four remaining coordinates
are terminal. Thus even in this example a coordinate indexed
by an upper entry need not equal that entry.

\subsection{Specialization and mixed quantum tori}

Write \(y_i\) for the normalized generic seed variables and \(Y_{r,c}\) for their specializations indexed by matrix positions. Let \(P\) be the integral unit upper triangular matrix whose \(j\)th column is the indicator of the chain up to \(j\). The leading exponent of \(y_j\) is \(Pe_j\), by triangular CGL construction. Let \(T_R\) be the initial mixed Laurent ring, which inverts mutable variables only. The integral quantum Laurent phenomenon, applied to the interval seeds, puts each \(x_i\) in \(T_R\), without requiring an a priori integral formula for all the CGL primes. The CGL triangular recursion expresses \(x_i\) with leading torus exponent \(P^{-1}e_i\) and Laurent-unit coefficient: the predecessor is mutable, and the correction involves only earlier generators. Here torus exponents are compared in reverse lexicographic order after applying \(P\). Products preserve this leading-order computation. Consequently
\[
 A_{\mathbb C(v)}\cap T_R=A_R.
\]
Indeed the highest PBW monomial of any element of the intersection gives a unique highest torus monomial with unit leading scalar, so its coefficient is in \(R\). Subtract it and continue over the finite PBW expansion. In particular every integral cluster variable, which lies in \(T_R\) and in \(A_{\mathbb C(v)}\), lies in \(A_R\); the converse inclusion follows from the interval seeds.

This saturation also identifies the kernel when \(T_R\) is specialized: an element of \(A_R\) divisible by \(v-\zeta^{(\ell+1)/2}\) in \(T_R\) is divisible by it already in \(A_R\). The same leading-term proof applies over \(\mathbb Z[v^{\pm1}]\):
the rewriting coefficients belong to this ring and the
normalizations and leading coefficients are powers of \(v\)
up to signs. It identifies the kernel of reduction by
\(\Phi_\ell(v)\), the \(\ell\)-th cyclotomic polynomial, before
extending scalars to \(\mathbb C\). The specialization and Laurent theorems of \citet[Theorems~3.10, 5.7 and 5.9]{NTY} therefore identify \(A_\zeta\) with the specialized quantum cluster algebra and embed it in the mixed Laurent ring of every seed. Compatibility diagonals have absolute value two, so there is no further condition on odd \(\ell\). In particular \(A_\zeta\) is a domain, by its embedding in the specialized torus.

With the compatibility diagonals of absolute value two,
the central power construction of Section~\ref{subsec:fiber-prelim}
applies in every specialized seed.
The original entries \(x_i\)
have central \(\ell\)-th powers because they occur in interval
seeds up to scalar. By the PBW basis, the algebra is finite
over the subalgebra generated by those powers, and hence also
over the full center. For a set \(S\) of frozen indices, let \(T_S\) denote the quantum polynomial--Laurent algebra in which precisely the coordinates in \(S\) are left uninverted. Then
\begin{equation}\label{eq:all-mixed-charts}
 A_\zeta[(Y_i^\ell)^{-1}:i\notin S]=T_S.
\end{equation}
The Laurent phenomenon gives containment from left to right. The other containment follows because the seed variables are actual elements and \(Y_i^{-1}=Y_i^{\ell-1}(Y_i^\ell)^{-1}\). The center also localizes as stated in \eqref{eq:all-mixed-charts}: one can clear the central denominators in a domain. Thus this is an equality of central localizations, not just an identification of fraction fields.

\subsection{The alternating radical and the boundary masks}\label{subsec:radical}

Index the pairs \((r,c)\), \(r<c\), in lexicographic order, and set
\(\beta_{(r,c)}=e_r+e_c\). The skew matrix of leading commutation
exponents has \(M_{ij}=(\beta_i,\beta_j)\) for \(i<j\), so that the
initial seed form is \(\Lambda=P^{\mathsf T}MP\).
Work over any commutative ring in which 2 is invertible. Given
\(v=(v_i)_{i=1}^N\), put
\[
 x_0=\tfrac12\sum_j v_j\beta_j,\qquad
 x_i=x_0-\sum_{j\le i}v_j\beta_j.
\]
Since \((\beta_i,\beta_i)=2\), we have
\begin{align*}
 (Mv)_i
 &=\sum_{j>i}v_j(\beta_i,\beta_j)-\sum_{j<i}v_j(\beta_i,\beta_j)\
 =2\bigl((\beta_i,x_{i-1})-v_i\bigr).
\end{align*}
Thus \(Mv=0\) implies \(x_i=s_{\beta_i}x_{i-1}\) and \(x_N=-x_0\).
Conversely, start with
\[
 x_N=s_{\beta_N}\cdots s_{\beta_1}x_0=-x_0,
 \qquad v_i=(\beta_i,x_{i-1}).
\]
Then \(\sum_i v_i\beta_i=x_0-x_N=2x_0\), and the same calculation
gives \(Mv=0\).

The reflection \(s_{e_i+e_j}\) replaces the pair of coordinates
\((a_i,a_j)\) by \((-a_j,-a_i)\). After the first \(k\) reflections
\((1,2),\ldots,(1,k+1)\), the first row has effected
\[
 (a_1,\ldots,a_n)\longmapsto
 (-a_{k+1},-a_1,a_2,\ldots,a_k,a_{k+2},\ldots,a_n).
\]
The remaining reflections act only on coordinates \(2,\ldots,n\).
Induction on \(n\), starting with \(n=2\), therefore gives
\[
 s_{\beta_N}\cdots s_{\beta_1}(a)=
 (-a_n,-a_{n-1},\ldots,-a_2,(-1)^{n+1}a_1).
\]
The equations for its anti-fixed vectors are
\[
 a_i=a_{n+1-i}\ (1\le i<n),\qquad a_n=(-1)^n a_1.
\]
For even \(n\) these are the mirrored-coordinate equations. For odd
\(n\) the first and last coordinates are zero and the others are
mirrored. In both cases the solution module is free of rank
\(\lfloor n/2\rfloor\).

Let \(F_r\) be the terminal variable of gap \(n-r\), with last pair
\((r,n)\), for \(1\le r\le n-2\). For \(n=2m\), let \(H,B\) denote
the terminal adjacency variables with leading matchings covering
\([1,n]\), \([2,n-1]\), respectively; for \(n=2m+1\), let \(B_L,B_R\)
denote those with matchings covering \([1,n-1]\), \([2,n]\).
The radical records exactly the exponent vectors giving central
torus monomials. Its basis below is supported on frozen coordinates
and separates into disjoint groups. If just one coordinate in a
group vanishes, a central monomial may kill that coordinate itself
in the fiber. If several vanish, the same central relation can
instead kill a product, leaving nonzero nilpotent classes.
These two possibilities give the factors in
Table~\ref{tab:chart-loewy}. The bases are:
\begin{align*}
 n=2m:\qquad & e_H,\quad e_B-e_{F_1},\quad
 e_{F_{n-r}}-e_{F_r}\quad(2\le r\le m-1);\\
 n=2m+1:\qquad &e_{B_L}+e_{B_R}-e_{F_1},\quad
 e_{F_{n-r}}-e_{F_r}\quad(2\le r\le m).
\end{align*}
To verify membership, fix \((i,j)\), write \(d=j-i\), and abbreviate
\(C_r=(MPe_{F_r})_{(i,j)}\). With \(g=n-r\), the terms from outgoing
edges at \(i,j\) and incoming edges at \(i,j\) give, respectively,
\begin{align*}
 C_r&=\operatorname{sgn}(g-d)\mathbf1_{i\le r}
      +\mathbf1_{j\le r}-\mathbf1_{i>g}
      -\operatorname{sgn}(g-d)\mathbf1_{j>g}.
\end{align*}
An edge identical to \((i,j)\) contributes zero, as reflected by
\(\operatorname{sgn}(0)=0\). Put \(s=n-r\), \(r\le s\), and
\[
 U=\mathbf1_{i\le r}-\mathbf1_{j>s},\qquad
 V=\mathbf1_{j\le r}-\mathbf1_{i>s}.
\]
By taking complements,
\(\mathbf1_{i\le s}-\mathbf1_{j>r}=V\) and
\(\mathbf1_{j\le s}-\mathbf1_{i>r}=U\). Therefore
\[
 C_r-C_s=(\operatorname{sgn}(s-d)-1)U+
          (1-\operatorname{sgn}(r-d))V=0.
\]
Indeed \(d\ge s\) implies \(i\le r\), \(j>s\), hence \(U=0\);
\(d\ge r\) implies \(j>r\), \(i\le s\), hence \(V=0\).
The same column formula at \(r=1\) yields
\(C_1=\mathbf1_{i=1}-\mathbf1_{j=n}\). For a matching of consecutive
pairs covering \(D\) the column is
\(-\mathbf1_{i\in D}+\mathbf1_{j\in D}\): if \(d>1\) the pair at
\(i\) precedes \((i,j)\), and the pair at \(j\) follows it.
If \(d=1\) and \((i,j)\) is itself a matching edge, the two
indicators cancel; otherwise the same signed rule applies. Consequently
the columns for \(H\), \(B-F_1\), and \(B_L+B_R-F_1\), in their
respective parities, are all zero.

Under the radical isomorphism above, in coordinates \(2x_0\),
the vector \(e_{F_{n-r}}-e_{F_r}\) maps to
\(2\mathbf1_{[r+1,n-r]}\), since the row-weight sum of a gap chain
is \(\mathbf1_{[1,r]}+\mathbf1_{[n-r+1,n]}\). In even size the
other two vectors map to \(\mathbf1_{[1,n]}\) and
\(\mathbf1_{[2,n-1]}-\mathbf1_{\{1,n\}}\); their sum and difference
give twice the interior and exterior layers. In odd size the triple
maps to \(2\mathbf1_{[2,n-1]}\). Subtracting successive nested
interval indicators now gives the individual mirrored layers. This
proves the basis assertion, in particular over \(\mathbb Z/\ell\mathbb
Z\). Mutation replaces \(\Lambda\) by \(E^{\mathsf T}\Lambda E\)
with \(E\) unimodular and \(Ee_F=e_F\) for every frozen index \(F\);
thus it preserves these bases in the seed coordinates.

The central exponent lattice of the full seed torus is the preimage
of this radical modulo \(\ell\). The lattice index, hence the rank
of the torus over its center, is \(\ell^{N-\lfloor n/2\rfloor}\).
The full seed torus is a central localization, so taking the square
root of this rank gives the PI degree (the degree of the central
simple quotient algebra)
\begin{equation}\label{eq:pidegree}
 \ell^{(N-\lfloor n/2\rfloor)/2}
 =\begin{cases}\ell^{m(m-1)}&n=2m,\\
               \ell^{m^2}&n=2m+1.
   \end{cases}
\end{equation}

\subsection{The full center}
\begin{proof}[Proof of Theorem~\ref{thm:intro-center}]
The radical calculation gives the asserted bases.
The PBW monomials \(x^{\ell w}\) are distinct, so
\(Z_0=\mathbb C[x_i^\ell:1\le i\le N]\) is polynomial. This is
exactly the Frobenius subalgebra: each \(x_i\) is a variable in an
interval seed up to a Laurent unit, and the classical counterparts
of the \(x_i\) generate the commutative algebra.

For ordered monomials \(x^a=x_1^{a_1}\cdots x_N^{a_N}\) put
\[
 \eta(a,b)=\sum_{j>i}a_jb_i M_{ji}.
\]
The leading straightening rule and its exponent difference read
\begin{align*}
 x^a x^b&=\zeta^{\eta(a,b)}x^{a+b}
               +\text{terms with exponent below }a+b,\\
 \eta(a,b)-\eta(b,a)&=a^{\mathsf T}Mb.
\end{align*}
The first identity follows by moving the \(b_i\) right-hand factors
across the \(a_j\) left-hand factors for \(j>i\); each correction
in the straightening has strictly smaller exponent. If
\(f=c x^a+\text{lower terms}\), \(c\ne0\), centrality gives
\[
 0=[x^{a+e_i}](fx_i-x_if)
  =c\zeta^{\eta(e_i,a)}(\zeta^{a^{\mathsf T}Me_i}-1)
 \quad(1\le i\le N).
\]
Thus \(Ma\equiv0\pmod\ell\). Since \(\Lambda=P^{\mathsf T}MP\),
this says \(k=P^{-1}a\bmod\ell\in K\), where \(K\) is the seed
radical modulo \(\ell\).

For each entry index \(i\), let \(f(i)\) be the terminal index
of its chain. The frozen columns of \(P\) are precisely the
disjoint chain indicators. Since \(\bar k\) is supported on
these columns,
\[
 (P\bar k)_i=\bar k_{f(i)},\qquad
 a_i\equiv\bar k_{f(i)}\pmod\ell,\qquad
 0\le\bar k_{f(i)}<\ell.
\]
Taking the quotient and remainder at each coordinate gives
\[
 a=P\bar k+\ell w,\qquad w\in\mathbb N^N .
\]
The based monomial \(Y^{\bar k}\) is in the algebra (frozen powers
are nonnegative) and is central in the torus.
If its leading coefficient is \(d_k\ne0\), the leading straightening
formula gives
\[
 [x^a](Y^{\bar k}x^{\ell w})
 =d_k\zeta^{\eta(P\bar k,\ell w)}=d_k .
\]
Thus replacing \(f\) by \(f-(c/d_k)Y^{\bar k}x^{\ell w}\)
kills its leading term and leaves a central polynomial with
strictly smaller leading exponent. The PBW order on
\(\mathbb N^N\) is well founded. The leading exponents
\(P\bar k+\ell w\) of the resulting family are all distinct,
since reduction modulo \(\ell\) followed by \(P^{-1}\) recovers
\(k\), and then subtraction and division recover \(w\).
This proves \eqref{eq:main-center}.
\end{proof}

For the center multiplication write
\[
 \bar k+\bar j=\overline{k+j}+\ell\sum_{F\ \mathrm{frozen}}
 \left\lfloor\frac{\bar k_F+\bar j_F}{\ell}\right\rfloor e_F.
\]
The based multiplication phase is one: the pairing of either
central representative with any exponent is divisible by \(\ell\),
and the chosen square root of \(\zeta\) has order \(\ell\).
In particular \(Y^{\ell e_F}=F^\ell\), and we obtain
\begin{equation}\label{eq:center-carry}
 Y^{\bar k}Y^{\bar j}=
 Y^{\overline{k+j}}
 \prod_{F\ \mathrm{frozen}}(F^\ell)^{
     \lfloor(\bar k_F+\bar j_F)/\ell\rfloor}.
\end{equation}

\section{Projective layers on seed charts}\label{sec:charts}

\begin{proof}[Proof of Theorem~\ref{thm:intro-chart}]
Use the notation of the theorem and let \(D\) be the quantum
torus on \(U\). Localization in \eqref{eq:all-mixed-charts}
does not change the fiber. With each \(Y^d\) ordered on \(S\), set
\[
 \begin{aligned}
 T_S&=\bigoplus_{d\in\mathbb N^S}D\,Y^d,\\
 \Gamma_0&=\{a\in\mathbb Z^U:a\bmod\ell\in K_0\},\\
 \Gamma_U&=\{a\in\mathbb Z^U:a\bmod\ell\in K_U\}.
 \end{aligned}
\]
Then \(C_0=\mathbb C[\Gamma_0]\) is the degree-zero center of
\(T_S\), and \(Z(D)=\mathbb C[\Gamma_U]\), in based
normalizations. For \(0\ne e\in\mathbb N^S\) with
\(e\bmod\ell\in\pi_S(K)\), choose a unit monomial \(u_e\) such
that \(u_eY^e\) is central. The full-center relations in
nonzero \(S\)-degree are determined by
\[
 \chi(u_eY^e)^\ell
     =\chi((u_eY^e)^\ell)=0,\qquad
 T_S\,u_eY^e=\bigoplus_{d\ge e} D\,Y^d .
\]
It follows that, with \(D_\chi=D\otimes_{C_0,\chi}\mathbb C\),
\begin{equation}\label{eq:chart-decomp}
 A_\chi=\bigoplus_{d\in\mathcal D_S}D_\chi\,Y^d .
\end{equation}

The inclusion \(\Gamma_0\subseteq\Gamma_U\) gives an
\'etale isogeny of tori (diagonalize the inclusion by Smith
normal form). Thus
\[
 Z(D)\otimes_{C_0,\chi}\mathbb C
  \cong\mathbb C^{\,|\Gamma_U/\Gamma_0|},\qquad
 |\Gamma_U/\Gamma_0|=|K_U|/|K_0|.
\]
Each factor of \(D_\chi\) is a twisted group algebra on
\(\mathbb Z^U/\Gamma_U\) with nondegenerate commutator.
Averaging over that finite group gives semisimplicity;
nondegeneracy gives one-dimensional center. Therefore
\[
 D_\chi\cong\prod_{\nu=1}^{|K_U|/|K_0|}M_s(\mathbb C),
 \qquad s^2=|\mathbb Z^U/\Gamma_U|
      =\ell^{|U|}/|K_U|.
\]

All \(d\in\mathcal D_S\) satisfy \(0\le d_i<\ell\), so the
positive-degree ideal in \eqref{eq:chart-decomp} is nilpotent.
Its quotient is \(D_\chi\). Using the fact that \(\mathcal D_S\)
is closed under decreasing coordinates, we obtain
\[
 J=\bigoplus_{|d|>0}D_\chi Y^d,\qquad
 J^b=\bigoplus_{|d|\ge b}D_\chi Y^d\quad(b\ge0),
\]
where the sums are over \(\mathcal D_S\). Indeed a degree
\(d\) monomial with \(|d|\ge b\) factors into \(b\) nonconstant
monomials for \(b\ge1\).
Write \(Y^d a=\sigma_d(a)Y^d\) for \(a\in D_\chi\).
The diagonal automorphism \(\sigma_d\) of the torus preserves
the reduction because it fixes \(C_0\). For a primitive
idempotent \(p_\alpha\) in \(D_\chi\), let \(\alpha(d)\) be the
simple label of \(D_\chi\sigma_d(p_\alpha)\). The nilpotent ideal
\(J\) preserves primitivity, and
\[
 P_\alpha=A_\chi p_\alpha,\qquad
 J^bP_\alpha/J^{b+1}P_\alpha
   \cong\bigoplus_{\substack{d\in\mathcal D_S\\|d|=b}}
              D_\chi\sigma_d(p_\alpha),\qquad
 C_{\beta\alpha}(t)
   =\sum_{\substack{d\in\mathcal D_S\\\alpha(d)=\beta}}t^{|d|}.
\]
In particular the sum over \(\beta\) is the asserted polynomial.

For the factorization, project the independent generators
of \(K\) onto \(S\). Disjoint supports make the standard
exponents a Cartesian product, with the following nontrivial
factors.
See Table~\ref{tab:chart-loewy}.
\begin{table}[htbp]
\caption{Standard exponents and radical-layer factors on a seed chart.}\label{tab:chart-loewy}
\centering
\begin{tabular}{>{$}c<{$}|>{$}c<{$}|>{$}l<{$}}
 \text{zero coordinates} & \text{standard exponents}
                      &\text{factor} \\ \hline
 \text{unpaired} & 0\le a<\ell & \sum_{a<\ell}t^a\\
 \text{opposite signs} & a+b<\ell & \sum_{j<\ell}(j+1)t^j\\
 \text{equal signs} & a,b<\ell,\ ab=0 &1+2\sum_{j=1}^{\ell-1}t^j\\[1pt]
 (-1,1,1)&
    a,b,c<\ell,\ a+\min(b,c)<\ell&
    \sum_{\text{indicated }a,b,c}t^{a+b+c}.
\end{tabular}
\end{table}
The exponents in the table are nonnegative. Each row counts the
monomials left on that support after the central relations are
imposed. If just one variable from a support vanishes, its
unit partners make the variable itself zero in the fiber, so
the factor is \(1\). For the three-zero
case the excluded minimal generators include the pure powers
and the vectors \((\ell-h,h,h)\), \(1\le h<\ell\); these give the
stated inequality. For the paired cases the vectors are
\((\ell-h,h)\) and \((h,h)\), respectively.

If a nontrivial factor occurs then \(J\ne0\), so \(A_\chi\)
is not a matrix algebra. Otherwise, for every \(i\in S\)
choose a unit monomial \(u_i\in D\) with \(u_iY_i\) central.
The central indeterminates \(v_i=u_iY_i\) give
\[
 D[v_i:i\in S]\xrightarrow{\ \sim\ }T_S,\qquad
 Z(T_S)=Z(D)[v_i:i\in S].
\]
Here injectivity follows by the \(S\)-degree decomposition.
The torus \(D\) is finite free over \(Z(D)\), with basis given
by the cosets of \(\Gamma_U\). Each fiber is a matrix
algebra by the same twisted-group-algebra argument, so \(D\)
and its polynomial extension are Azumaya over their centers.
The number of admissible masks is
\(2\cdot 3^{m-1}\) in even size (the singleton \(H\) and
\(m-1\) pairs), and \(4\cdot3^{m-1}\) in odd size
(the triple and \(m-1\) pairs).
\end{proof}

\begin{example}
For \(n=4\) and \(\ell=3\), the opposite-sign pair is
\((F_1,B)\), and \(F_2\) is unpaired. The two zero sets give
\[
 \{F_1,B\}:\ 1+2t+3t^2,
 \qquad
 \{F_2\}:\ 1+t+t^2 .
\]
The same Loewy length \(3\) corresponds to different
composition lengths. If only \(F_1\) vanishes, the fiber
is a matrix algebra.
\end{example}

\subsection{The global Azumaya locus}
Let \(Z=Z(A_\zeta)\), let \(\operatorname{Reg}(Z)\) be the set
of primes with regular local ring, and write \(D(g)\) for the
principal open on which \(g\) is a unit. We call a prime Azumaya
if the central localization there is an Azumaya algebra.

\begin{corollary}\label{cor:global-azumaya}
For all the odd orders considered here,
\[
 \operatorname{Az}(A_\zeta)=
 \begin{cases}
 \operatorname{Reg}(Z),& n\text{ odd},\\
 \operatorname{Reg}(Z)\cap D(F_m^\ell),&n=2m.
 \end{cases}
\]
Here \(F_m=Y_{m,n}\) is the unpaired frozen variable in even size.
\end{corollary}

\begin{proof}
Put \(P=A_R/(v-1)A_R\). The integral PBW basis identifies \(P\)
with the polynomial algebra on the upper entries. In an interval
presentation the semiclassical bracket has the form
\[
 \begin{aligned}
  \{\bar x_j,\bar x_i\}
   &=\lambda_{ji}\bar x_i\bar x_j+\partial_j(\bar x_i),\\
  \partial_j(\bar x_i)&=2\sum_{s,t}c_{st}\bar x_s\bar x_t,\\
  \partial_j(\bar x_s)&=\partial_j(\bar x_t)=0
     \quad(c_{st}\ne0).
 \end{aligned}
\]
Here the coefficients are the specializations of the quantum
corrections after removing \(q-q^{-1}\), since
\(\left.(q-q^{-1})/(q-1)\right|_{q=1}=2\).
The two factors in each correction share an endpoint with the added
generator. Thus \(\partial_j^2(\bar x_i)=0\) on generators, and the
Leibniz rule gives local nilpotence. The forward and reverse torus
weights at the added generator are \(-2\) and \(2\), respectively.
This gives the symmetric Poisson--CGL presentation. Distinct
ordered products remain distinct at \(v=1\), so the specialized
derivations have the same zero set of chain starts.

For an integral interval prime \(y\), the monic PBW term gives
\(\bar y\ne0\). Write \(hy=yh'\) in the generic algebra. With
\(\nu\) the coefficient valuation at \(v=1\), multiplication by
\(y\) preserves \(\nu\), by the domain property of \(P\). Hence
\[
 \begin{gathered}
  \nu(h')=\nu(yh')=\nu(hy)\ge0,\qquad
  \bar y(\overline{h'}-\bar h)=0,\\
  \{\bar h,\bar y\}
   =\bar y\,\overline{(h'-h)/(q-1)}.
 \end{gathered}
\]
The triangular correction in the monic recursion is integral and
retains its earlier PBW support. Thus \citet[Proposition~5.7]{GYPoisson}
identifies the specialized sequences as the Poisson prime sequences
with the same mutable positions. The power-seed map sends each
interval variable \(\bar y\) to \(Y^\ell\) by the exchange rule.
In particular it sends \(\bar x_{ij}\) to \(x_{ij}^\ell\):
the signs of the normalizing units \(\pm v^j\) agree at \(v=1\)
and after taking the odd \(\ell\)-th power. This identifies
\(P\cong Z_0\) and the corresponding interval opens.

Let \(d_\tau\) be the product of mutable powers in an interval
chart. By \citet[Theorem~10.5(d)]{GYPoisson},
\[
 Z_0=\bigcap_{\tau\in\Gamma_N}(Z_0)_{d_\tau},\qquad
 d_\tau\in(p)\ \text{for all }\tau
 \ \Longrightarrow\
 p^{-1}\in\bigcap_{\tau\in\Gamma_N}(Z_0)_{d_\tau}\setminus Z_0
\]
for any height-one prime \((p)\subset Z_0\); the implication is a
contradiction. For \(\mathfrak p\subset Z\) of height one,
incomparability and going-down for the finite extension of the
normal domain \(Z_0\) give
\(\operatorname{ht}(\mathfrak p\cap Z_0)=1\).
Hence \(\mathfrak p\) is in an interval chart. In that chart
the frozen powers are independent variables, so at most one can
belong to \(\mathfrak p\). Let \(S(\mathfrak p)\) be the corresponding frozen zero set.
On a dense open of \(V(\mathfrak p)\) in the chart the frozen zero
set of every closed point is \(S(\mathfrak p)\). By
Theorem~\ref{thm:intro-chart} and openness,
\[
 |S(\mathfrak p)|\le1,\qquad
 \bigl(n\text{ odd or }F_m^\ell\notin\mathfrak p\bigr)
 \quad\Longrightarrow\quad \mathfrak p\in\operatorname{Az}(A_\zeta).
\]

The automorphic Ore presentation gives
\(\operatorname{gl.dim}A_\zeta\le N\), also after central
localization. For each maximal ideal \(\mathfrak m\subset Z\),
let \(z_1,\ldots,z_N\) be regular parameters of
\((Z_0)_{\mathfrak m\cap Z_0}\). Freeness over \(Z_0\) and
localization show that each multiplication map
\[
 z_i:A_{\mathfrak m}/(z_1,\ldots,z_{i-1})A_{\mathfrak m}
       \longrightarrow A_{\mathfrak m}/(z_1,\ldots,z_{i-1})A_{\mathfrak m}
\]
is injective. The last quotient is nonzero by Nakayama, so
\(\operatorname{depth}_{Z_{\mathfrak m}} A_{\mathfrak m}
=N=\dim Z_{\mathfrak m}\).
This verifies the Cohen--Macaulay hypothesis of
\citet[Theorem~3.13]{BM}, which gives equality of the two loci
in odd size and on \(D(F_m^\ell)\) in even size. In the even case,
each component of \(V(F_m^\ell)\) has height one and its generic
point is non-Azumaya by the single unpaired zero case of the
chart criterion. Since the non-Azumaya locus is closed, it
contains \(V(F_m^\ell)\). This proves the asserted equality.
\end{proof}

\section{Exchange localizations and boundary fibers}
\label{sec:endpoint}

\subsection{Unit exchange monomials}

\label{subsec:unit-exchange-patch}

Fix a seed of $A_\zeta(n)$ and a mutable direction, with variable $X$
and mutated variable $W$. Let $A_1,A_2$ be the two based exchange
monomials in the product $XW=A_1+A_2$. Invert the central
$\ell$th powers of every shared mutable variable and of the frozen
variables occurring in $A_1$ or $A_2$. Further spectator frozen powers
may be inverted or left uninverted. Denote the resulting localization by
$\mathcal L$ and its shared mixed quantum torus by $D$.

\begin{theorem}\label{thm:unit-exchange-patch}
Suppose the compatibility diagonal at the chosen direction is coprime
to $\ell$. Put $t=A_1A_2^{-1}$, $\sigma=\operatorname{Ad}_X|_D$,
$D_0=D^\sigma$, $C=Z(D)^\sigma$, and $\rho=t^\ell$.
Then $Z(D_0)=C$, and $\mathcal L\simeq D_0\otimes_C Q$, where $Q$
is generated over $C$ by $t,X,Z$, with $Z=WA_2^{-1}$ and
\[
 \begin{gathered}
 t^\ell=\rho,\qquad Xt=\omega tX,\qquad Zt=\omega^{-1}tZ,\\
 XZ=1+t,\qquad ZX=1+\omega^{-1}t,
 \end{gathered}
\]
for a primitive $\ell$th root $\omega$. The full center of $\mathcal L$ is
\[
 S=C[\mathsf x,\mathsf z]/(\mathsf x\mathsf z-1-\rho),
 \qquad \mathsf x=X^\ell,\quad \mathsf z=Z^\ell .
\]
The algebra $Q$ is Azumaya of degree $\ell$ over $S$.
For each closed point $m$ of $S$, with $c=m\cap C$, the full-center
fiber is $M_\ell(D_0\otimes_C\kappa(c))$. Consequently it is
isomorphic, as an algebra, to the fiber at
$(c,\mathsf x=1,\mathsf z=1+\rho(c))$ in the $X$-inverted seed
chart. The Azumaya locus of $\mathcal L$ is the inverse image of
that of $D_0$.
\end{theorem}

\begin{proof}
The difference of the exchange exponents is, up to sign, the exchange
column. Its skew pairing with each shared direction is zero and its
pairing with $X$ is the compatibility diagonal up to sign. Thus
$t\in Z(D)^\times$ and $\sigma(t)=\omega t$. On the shared
generators $\sigma$ is diagonal with $\ell$th-root eigenvalues, so
$\sigma^\ell=1$. In the $X$-Laurent ring we have
$Z=X^{-1}(1+t)$, giving the displayed relations and
$Zd=\sigma^{-1}(d)Z$ for $d\in D$.

By \eqref{eq:all-mixed-charts}, $\mathcal L$ lies in the
intersection of the two mixed Laurent rings on $X$ and $Z$; replacing
$W$ by $Z$ does not change the second ring. The signed decomposition
can be seen directly, or compared with the rank-one calculation of
\citet[Theorem~4.2]{HLY}. In the common fraction ring,
\[
 Z^j=\left(\prod_{a=1}^j(1+\omega^{-a}t)\right)X^{-j}
 \qquad(j>0).
\]
On comparing degrees in $X$, negative-degree coefficients in the
intersection give exactly $DZ^j$; positive-degree coefficients give
$DX^j$, since $X^jZ^j\in D$. All these signed summands are already
in $\mathcal L$. Hence
\[
 \mathcal L=D\oplus\bigoplus_{j>0}(DX^j\oplus DZ^j).
\]
The directness follows from the distinct Laurent degrees and the
domain property. Multiplication by $t^i$ identifies $D_0$ with
the $\omega^i$-eigenspace, giving
$D=\bigoplus_{i=0}^{\ell-1}D_0t^i$. Since $t$ is central in $D$,
an element of $Z(D_0)$ commutes with all of $D$ and is fixed by
$\sigma$. The converse is immediate. The generators of $Q$ commute
with $D_0$, and the signed decomposition, using the free basis
$1,t,\ldots,t^{\ell-1}$ of $C[t]/(t^\ell-\rho)$, proves the
tensor identification.

Commuting with $t$ restricts a central element to signed degrees
divisible by $\ell$. The generators $X^\ell$ and $Z^\ell$ in
these degrees are central: $\sigma^\ell=1$, and the product along
a full orbit is
\[
 X^\ell Z^\ell=\prod_{i=0}^{\ell-1}(1+\omega^it)=1+\rho .
\]
Here $\ell$ is odd. For a coefficient $d$ in such a signed component,
commuting with $D$ gives $d\in Z(D)$, and commuting with $X$ gives
$\sigma(d)=d$, by cancellation in the mixed domain. Thus the
signed basis proves the assertion about the whole center, with no
extra relation. The same calculation gives $Z(Q)=S$.

Adjoin $\lambda^\ell=\rho$ to $S$, a faithfully flat finite
\'etale extension since $\rho$ is a unit. The $t$-eigenspaces
split into $\ell$ orthogonal idempotents with eigenvalues
$\omega^i\lambda$. At any geometric point at most one factor
$1+\omega^i\lambda$ vanishes. On the open set where the other
$\ell-1$ factors are invertible, the corresponding bonds between
successive idempotents are invertible, with inverses given by the
reverse arrows divided by these factors. Those bonds form a connected
path and yield matrix units. Each diagonal corner equals the extended
base ring: a signed degree not divisible by $\ell$ moves to a
different idempotent, while a signed degree divisible by $\ell$ is
a central power times a degree-zero term. The degree-zero corner of
$C[t]/(t^\ell-\rho)$ after base extension is scalar. Thus $Q$ is
\'etale-locally a full matrix algebra. This proof includes the locus
$\mathsf x=\mathsf z=0$.

The resulting matrix description gives the asserted closed fibers.
The indicated point with $\mathsf x=1$ lies in the $X$-inverted chart,
and both points have coefficient character $c$. Finally $D_0$ is
finite over $C$ by the shared Frobenius powers and the direct-summand
decomposition of $D$. The map $C\to S$ is faithfully flat, by the
basis $1,\mathsf x^j,\mathsf z^j$ for $j>0$. After the
\'etale matrix trivialization of $Q$, invariance under passing to
a matrix algebra and faithfully flat descent give the stated
Azumaya locus.
\end{proof}

The compatibility diagonals have absolute value two in the skew-matrix
atlas, so the theorem applies at every odd order. It compares fibers;
chart coverage of a fixed point is a separate question. Since $A_2$
quasi-commutes with $W$, the powers $Z^\ell$ and $W^\ell$ differ
by the central unit $A_2^{-\ell}$ in the odd-order convention.

\subsection{Seed coverage with one mutable coordinate allowed to vanish}
\label{subsec:single-mutable-coverage}

In this subsection all exchanges are taken in the commutative Frobenius
algebra. We write \(y_{ij}=Y_{i,j}^{\ell}\) and use the same letters
for their values when evaluating at a point. For a mutable index
\((i,j)\), let \(U_{ij}\) be the open where all other initial mutable
powers are nonzero. No frozen power is inverted in the definition of
\(U_{ij}\).

\begin{proposition}\label{prop:single-mutable-coverage}
Suppose \(n\geq6\), and write \(yv=A+B\) for the power exchange at
\(y=y_{ij}\). On \(U_{ij}\), the initial or once-mutated seed covers
every point with \(y\ne0\) or \(v\ne0\). At a point with \(y=v=0\),
the following alternatives hold.
\begin{enumerate}
\item[(i)] If \(A,B\ne0\), a mixed seed reached by two mutations
covers the point.
\item[(ii)] If \(A=B=0\), no mixed seed covers the point. In the
fixed-frozen Frobenius fiber its local model is a torus times the
node \(yv=0\). The zero conditions \(A=B=0\) on \(U_{ij}\)
are precisely the following:
\[
\begin{array}{c|l}
 (i,j) & \text{frozen powers at the point}\\ \hline
 j=n-1,\ j-i\geq3 & y_{i,n}=y_{i+1,n}=0\\
 (i,j)=(n-3,n-1) & y_{n-2,n}=0,\quad y_{n-3,n}y_{n-2,n-1}=0\\
 (i,j)=(n-3,n-2) & y_{n-2,n}=y_{n-1,n}=0 .
\end{array}
\]
\end{enumerate}
The same coverage criterion holds for full-center characters above
these points.
\end{proposition}

\begin{proof}
First compute the exchange columns. Let \(\pi(i,j)\) denote the
chain predecessor of \((i,j)\), equal to \((i-1,j-1)\) for
\(j-i\geq2\) and \((i-2,i-1)\) for \(j-i=1\); nonexistent
predecessors contribute zero to exponent vectors. The
successor-negative column \(b_{ij}\) satisfies
\begin{equation}\label{eq:all-exchange-columns}
P b_{ij}=
\begin{cases}
 e_{(i,j+1)}+e_{(i+1,j)}-e_{(i,j)}-e_{(i+1,j+1)},& j-i\geq2,\\
 -e_{(i,i+1)}+e_{(i,i+2)}+e_{(i+1,i+3)}-e_{(i+2,i+3)},& j-i=1.
\end{cases}
\end{equation}
To verify this, write \(u\) for either right-hand side,
\(\beta_{ab}=\epsilon_a+\epsilon_b\), and use the lexicographic order
on positions. Both candidates have \(\sum_s u_s\beta_s=0\).
At an index \(t\) outside the support, the shared-endpoint formula
for \(M\) gives
\[
 u^{\mathsf T}Me_t=2\Bigl(\sum_{s<t}u_s\beta_s,\beta_t\Bigr).
\]
For gap at least two, the only interval with intervening positions
and nonzero partial weight has partial weight
\(\epsilon_{j+1}-\epsilon_j\), between \((i,j+1)\) and
\((i+1,j)\). For gap one, the corresponding intervals are
between \((i,i+2)\) and \((i+1,i+3)\), and between
\((i+1,i+3)\) and \((i+2,i+3)\), with partial weights
\(\epsilon_{i+2}-\epsilon_{i+1}\) and
\(\epsilon_{i+2}+\epsilon_{i+3}\), respectively. In each interval
the intervening roots are orthogonal to the displayed weight.
On the four support indices, in increasing order, the pairings
are \((-2,0,0,2)\). The first and fourth indices are the
current coordinate and its successor on the same chain. Thus \(u^{\mathsf T}MP=-2e_{(i,j)}^{\mathsf T}\).
Compatibility and weight zero determine the column: the
difference of two candidates, after multiplication by \(P\), lies
in \(\ker M\) and has zero weight, hence is zero by the weight
injectivity on the commutation radical.

Since \(P^{-1}e_{ab}=e_{ab}-e_{\pi(a,b)}\), the resulting
columns are
\[
\begin{aligned}
 b_{ij}={}&e_{(i,j+1)}+e_{(i+1,j)}+e_{\pi(i,j)}\\
          &-e_{(i+1,j+1)}-e_{\pi(i,j+1)}-e_{\pi(i+1,j)}
              \qquad (j-i\geq2),\\
 b_{i,i+1}={}&e_{(i+1,i+3)}+e_{\pi(i,i+1)}
       -e_{(i+2,i+3)}-e_{\pi(i,i+2)}.
\end{aligned}
\]
The column computation itself also holds for the mutable
indices when \(n=5\). In particular, the first-row columns are
\begin{equation}\label{eq:first-row-column}
 b_{1j}=e_{(1,j+1)}+e_{(2,j)}-e_{(2,j+1)}
       -\mathbf1_{\{j>3\}}e_{(1,j-1)}\qquad (3\leq j\leq n-1).
\end{equation}
Among the factors in these columns, the only frozen factors are
the ones just listed in the statement and the negative gap-one
factor \(y_{n-2,n-1}\) in the exchange at \((n-4,n-3)\). In
the last case the positive exchange monomial is a unit, so both
terms cannot vanish. This proves the list in (ii).

For (i), the following choices of a shared mutable \(z=y_{ab}\)
suffice:
\[
\begin{array}{l|c}
 (i,j) & (a,b)\\ \hline
 j-i\geq2,\ i\geq2,\ j\leq n-2 &(i-1,j)\\
 i=1,\ 4\leq j\leq n-2 &(1,j-1)\\
 (i,j)=(1,3)&(1,4)\\
 j=n-1,\ j-i\geq3 &(i,j-1)\\
 (i,j)=(n-3,n-1)&(i-1,i)\\
 j=i+1,\ i\leq n-5 &(i+1,i+3)\\
 j=i+1,\ i\in\{n-4,n-3\}&(i-1,i+1).
\end{array}
\]
The columns above give \(b_{z,y}=-b_{y,z}\in\{1,-1\}\) in
each row. After removing \(y\), the exchange at \(z\) has only
mutable factors except for \(y_{n-2,n-1}\) in the fifth row.
There this frozen factor also occurs in the exchange at \(y\)
and is nonzero by hypothesis. Hence, orienting the signs so that
the monomial containing \(z\) is first, the exchanges read
\[
 yv=za+b,\qquad zz'=c+yd
\]
with \(a,b,c,d\) nonzero at the point. Column mutation at \(z\)
gives \(yu=(ac+bz')/g\) for the next exchange at \(y\), where
\(g=\gcd(c,b)\) is the common monomial factor removed from
the two terms. Thus
\[
 zyu=(c/g)(za+b)+y(b/g)d,\qquad
 zu=(c/g)v+(b/g)d .
\]
Cancellation takes place in the Frobenius domain before
evaluation. At \(y=v=0\) the new powers are
\(z'=c/z\) and \(u=(b/g)d/z\), both nonzero. All other
mutable powers are unchanged, so this seed covers the point.

Finally let \(R\) be the Laurent polynomial ring in the shared
mutable powers with polynomial frozen coefficients, and put
\(p=A+B\). In the fraction field of the power algebra one has
\[
 R[y^{\pm1}]\cap R[v^{\pm1}]=R[y,v]/(yv-p),\qquad v=p/y.
\]
Indeed, a negative \(y\)-degree \(-d\) in the intersection
has coefficient divisible by \(p^d\); all nonnegative terms
already lie in \(R[y]\). The Laurent inclusions in both adjacent
seeds and the presence of \(R,y,v\) give equality with the actual
partial localization. Fixing the frozen values when \(A=B=0\)
now gives the claimed torus times a node, singular at \(y=v=0\).
In any mixed seed the fixed-frozen fiber is smooth, since its
coordinate ring is Laurent polynomial in the mutable powers.
Thus no such seed can contain this point. Membership in the
mixed opens is decided by the powers and therefore is unchanged
on passage to a full-center character.
\end{proof}

\subsection{The terminal exchange and the nodal fiber}

In this section let \(n\geq5\), and let the specialization parameter be a primitive root of unity of odd order \(\ell>1\) over \(\mathbb C\). We use the initial seed and its mixed Laurent charts established above. Write \(F_r=Y_{r,n}\) for the terminal variable of the gap \(n-r\), \(1\leq r\leq n-2\). Localize at the central powers of all initial mutable variables except \(X=Y_{1,n-1}\), and denote its first mutation by \(W\). At a fixed frozen zero set \(S\), also invert the powers of the frozens outside \(S\). Write \(D\) for the resulting shared mixed quantum torus on the variables other than \(X,W\), and write \(\mathsf{x}=X^\ell\), \(\mathsf{w}=W^\ell\) for
central elements. We write \(\chi(g)\) for their actual values; in
fiber equations, \(g=0\) abbreviates \(\chi(g)=0\). The exchange has the form
\begin{equation}
 XW=A+B,\qquad WX=\sigma^{-1}(A+B),
 \label{eq:endpoint-quantum}
\end{equation}
where \(\sigma\) is conjugation by \(X\) in the quotient division ring. The monomial \(A\) is supported on \(F_1,Y_{2,n-1}\), and \(B\) on \(F_2,Y_{1,n-2}\), with scalar normalizations as in the based torus. In particular the commutative Frobenius exchange is \(\mathsf{x}\mathsf{w}=af_1+yf_2\), where \(a,y\) are the powers of \(Y_{2,n-1},Y_{1,n-2}\), respectively. A \emph{collision} means \(\mathsf{x}=\mathsf{w}=0\).

\begin{theorem}
\label{thm:endpoint}
In the endpoint partial localization the following statements hold.
\begin{enumerate}
\item[(i)] If exactly one of \(f_1,f_2\) vanishes, a collision is impossible. If both are nonzero, every collision belongs to the seed obtained by mutating first at \(Y_{1,n-2}\) and then at \(X\). Such a fiber is Azumaya exactly when its unchanged frozen zero set is an Azumaya mask for a seed chart.
\item[(ii)] If \(f_1=f_2=0\), the collision lies outside every seed chart. Such central characters exist for every choice of the spectator zero mask and every character of the degree-zero coefficient center described below.
\item[(iii)] Suppose \(f_1=f_2=0\). The partners of \(F_1\)
in its radical support are the frozen terminal \(B\) from
Section~\ref{sec:intro} in even size and \(B_L,B_R\)
in odd size. The partner of \(F_2\) is \(F_{n-2}\).
Let \(\mathbf{s}\) be a multi-index on the other zero frozen
coordinates. Use the ordinary mixed-chart monomial bounds in every radical group disjoint from \(F_1,F_2\). In the two groups containing \(F_1,F_2\), let \(p=p(\mathbf{s})\), \(r=r(\mathbf{s})\) be the numbers of allowed exponents of \(F_1,F_2\), respectively, when the spectator exponents are fixed. Explicitly,
\begin{equation}
\begin{aligned}
p&=1 &&\text{if no partner of }F_1\text{ is zero},\\
p&=\ell-s_G &&\text{if exactly one partner }G\text{ is zero},\\
p&=\ell-\min(s_{G_1},s_{G_2}) &&\text{if both partners are zero (odd }n\text{)},\\
r&=1 &&\text{if }F_{n-2}\text{ is a unit},\\
r&=\ell-s_{F_{n-2}} &&\text{if }F_{n-2}\text{ is zero}.
\end{aligned}
\label{eq:endpoint-bounds}
\end{equation}
All spectator exponents are less than \(\ell\). For \(|k|<\ell\) and \(j\geq0\), let \(m=\ell-|k|\) and set
\begin{equation}
\begin{aligned}
 R_j &= ([0,j-p]\cup[r,j])\bmod\ell,\\
 h_j(\mathbf{s},k)&=\min(m,j+1)-\min(m,|R_j|),
\end{aligned}
\label{eq:endpoint-rank}
\end{equation}
where intervals are integer intervals and are empty if the upper endpoint is below the lower endpoint.
Every projective cover in the fiber has the same unlabelled Loewy polynomial
\begin{equation}
 \sum_{\mathbf{s}}\sum_{|k|<\ell}\sum_{j\geq0}
 h_j(\mathbf{s},k) z^{|\mathbf{s}|+|k|+2j}.
 \label{eq:endpoint-loewy}
\end{equation}
Its coefficient at \(z^d\) is the composition length of the \(d\)-th radical layer. In each summand, the simple label is obtained by applying the common torus-character permutation of the corresponding monomials.
\end{enumerate}
\end{theorem}

\begin{proof}[Proof of (i) and (ii)]
Use the first-row columns \eqref{eq:first-row-column} at
\(c=n-2,n-1\). They give the support in
\eqref{eq:endpoint-quantum}.

Let \(d=Y_{2,n-2}^{\ell}\), and set \(e=Y_{1,n-3}^{\ell}\) for
\(n\geq6\), \(e=1\) for \(n=5\), using these symbols also for their
values. The elements \(a,y,d,e\) are units in the chosen
localization. The
initial exchange columns have signs \(f_1+a-f_2-y\)
and \(\mathsf{x}+d-a-e\), with symbols here denoting coordinate
indices (omit \(e\) when \(n=5\)). Mutation at \(y\) changes
the former by replacing the row \(-y\) by \(+v\) and
subtracting \(a+e\) on the other rows. Thus the column
at \(\mathsf{x}\) changes to \(f_1+v-f_2-e\).
The corresponding power exchanges are
\[
 \mathsf{x}\mathsf{w}=af_1+yf_2,\qquad yv=\mathsf{x}d+ea,\qquad
 \mathsf{x}u'=vf_1+ef_2,
\]
where \(v\) replaces \(y\) and \(u'\) then replaces \(\mathsf{x}\).
Eliminate the unwanted term before evaluating at the collision:
\[
 \begin{aligned}
 a\mathsf{x}u' &= v(af_1)+aef_2\\
      &= v(\mathsf{x}\mathsf{w}-yf_2)+aef_2\\
      &= \mathsf{x}\mathsf{w}v+(ae-yv)f_2
       =\mathsf{x}(\mathsf{w}v-df_2).
 \end{aligned}
\]
The commutative Frobenius algebra is a domain, so
\(au'=\mathsf{w}v-df_2\). At \(\mathsf{x}=\mathsf{w}=0\) the necessary and resulting
identities are
\[
 af_1+yf_2=0,\qquad v=ea/y,\qquad u'=-df_2/a.
\]
As \(a,y\) are units, exactly one frozen zero cannot satisfy the
first equation. If both frozen values are units, then \(v,u'\)
are units and the collision lies in the indicated double mutation
chart. The frozen mask is unchanged, so the previous chart
criterion applies there, proving (i).

For (ii), before specialization write \(R\) for the commutative
coefficient ring of the partial localization (with the shared
mutable variables inverted), and set \(p=af_1+yf_2\).
The two Laurent charts intersect in
\[
 R[\mathsf{x}^{\pm1}]\cap R[(p/\mathsf{x})^{\pm1}]
   =\bigoplus_{j\geq0}R\mathsf{x}^j\ \oplus\
       \bigoplus_{j>0}R(p/\mathsf{x})^j
   =R[\mathsf{x},\mathsf{w}]/(\mathsf{x}\mathsf{w}-p).
\]
Indeed a term of \(\mathsf{x}\)-degree \(-j<0\) in the second chart has
coefficient divisible by \(p^j\), while a term of nonnegative
degree already lying in the first chart lies in the second.
The Laurent property bounds the actual Frobenius localization
above by this intersection, and it contains \(R,\mathsf{x},\mathsf{w}\).
After fixing the frozen values with \(f_1=f_2=0\), this gives
\[
 \begin{gathered}
 \mathbb C[z_1^{\pm1},\dots,z_h^{\pm1},\mathsf{x},\mathsf{w}]/(\mathsf{x}\mathsf{w}),\\
 \dim T_PX_c=h+2>h+1=\dim_P X_c
 \quad(\mathsf{x}(P)=\mathsf{w}(P)=0).
 \end{gathered}
\]
where \(X_c\) here denotes the fixed-frozen Frobenius fiber and the
\(z_i\) are its other initial mutable coordinates. In a full seed
chart that same fixed-frozen fiber is a smooth torus, so the
collision is in no seed chart. The chart noncoverage also holds
on the full center, since every mutable unit condition is detected
by its central power. Existence over the specified degree-zero
characters follows in the signed calculation below.
\end{proof}

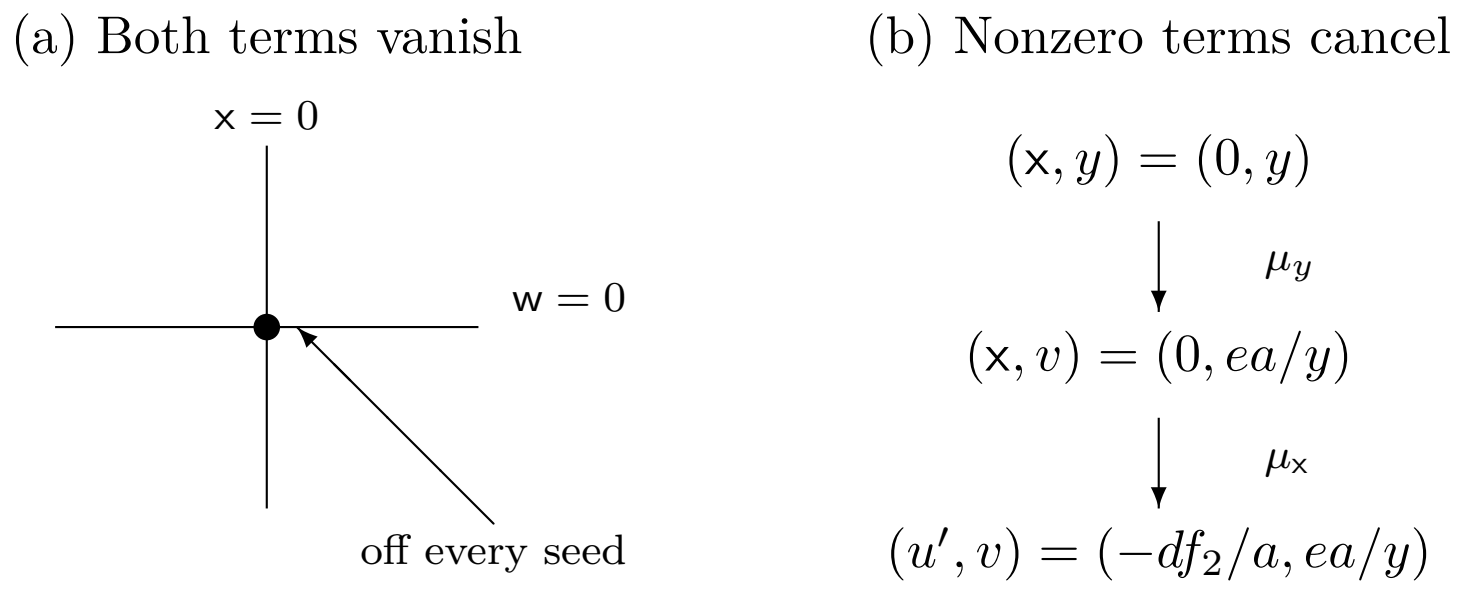
\begin{figure}[htbp]
\centering
\setlength{\unitlength}{1mm}
\begin{picture}(119,43)
\put(28,40){\makebox(0,0){\small (a) Both terms vanish}}
\put(87,40){\makebox(0,0){\small (b) Nonzero terms cancel}}
\put(14,21){\line(1,0){28}}
\put(28,9){\line(0,1){24}}
\put(28,21){\circle*{1.8}}
\put(28,35){\makebox(0,0){\scriptsize \(\mathsf{x}=0\)}}
\put(48,23){\makebox(0,0){\scriptsize \(\mathsf{w}=0\)}}
\put(43,8){\vector(-1,1){13}}
\put(43,6){\makebox(0,0){\scriptsize off every seed}}
\put(87,32){\makebox(0,0){\small \((\mathsf{x},y)=(0,y)\)}}
\put(87,28){\vector(0,-1){6}}
\put(94,25){\makebox(0,0)[l]{\scriptsize \(\mu_y\)}}
\put(87,19){\makebox(0,0){\small \((\mathsf{x},v)=(0,ea/y)\)}}
\put(87,15){\vector(0,-1){6}}
\put(94,12){\makebox(0,0)[l]{\scriptsize \(\mu_{\mathsf{x}}\)}}
\put(87,6){\makebox(0,0){\small \((u',v)=(-df_2/a,ea/y)\)}}
\end{picture}
\caption{Terminal collision: two zero exchange terms give a fixed-frozen
node (left); cancellation of two units gives a chart after two mutations (right).}
\label{fig:endpoint}
\end{figure}

Figure~\ref{fig:endpoint} contrasts the fixed-frozen node with the two-mutation route.

\begin{proof}[Proof of (iii)]
Set \(U=A+B\), and work in the common quotient division ring.
Write \(T\) for the endpoint partial localization and put
\[
 U_j=\sigma^j(U),\qquad
 N_t^+=\prod_{j=0}^{t-1}U_j,\qquad
 N_t^-=\prod_{j=1}^{t}U_{-j}.
\]
The exponent difference of \(A,B\) pairs trivially with
the shared variables and by \(\pm2\) with \(X\).
Thus \(A,B\) commute, their conjugations on \(D\) agree,
and each \(U_j\) is a normal element of \(D\).
Moreover \(\sigma'(d)=U_{-1}\sigma^{-1}(d)U_{-1}^{-1}\)
for conjugation by \(W\), so
\(\sigma'(U_j)=U_{j-1}\). Consequently
\[
 W=U_{-1}X^{-1},\qquad
 W^t=N_t^-X^{-t},\qquad
 X^tW^t=N_t^+,\quad W^tX^t=N_t^- .
\]
Comparing homogeneous coefficients in the two mixed Laurent
charts gives, for \(t>0\),
\[
 \begin{aligned}
 DX^t&\subseteq DW^{-t}           &&(X^t=N_t^+W^{-t}),\\
 DW^t&\subseteq DX^{-t}           &&(W^t=N_t^-X^{-t}),\\
 \bigl(D[X^{\pm1};\sigma]\cap D[W^{\pm1};\sigma']\bigr)_{\!k}
   &=
 \begin{cases}
   DX^k,& k\geq0,\\
   DW^{-k},& k<0 .
 \end{cases}
 \end{aligned}
\]
The Laurent property gives containment of \(T\) in the intersection,
and \(T\) contains \(D,X,W\). Hence \(T=\bigoplus_k T_k\) with
these summands, and this decomposition survives base change over
the degree-zero center.

Write \(\sigma(A)=\alpha A\), \(\sigma(B)=\beta B\) and
\(\omega=\alpha/\beta\). The scalars \(\alpha,\beta\) are \(\ell\)-th
roots of unity, and \(\omega\) has order \(\ell\).
Since \(A,B\) commute, \(U_j=\beta^j(B+\omega^j A)\);
oddness of \(\ell\) gives
\[
 N_\ell^+=N_\ell^-=A^\ell+B^\ell .
\]
This also follows up to a nonzero scalar from the set of
\(\ell\) relative phases, which suffices for the central-power
ideal below. With \(t=AB^{-1}\) in the division ring,
\[
 [t,D]=0,\qquad
 X t=\omega t X,\qquad Wt=\omega^{-1}tW.
\]
A central element commutes with \(t\); on a homogeneous element of
signed degree \(k\) the commutation factor is \(\omega^k\).
It follows that central degrees are multiples of \(\ell\).
In those degrees \(X^{\ell}\) and \(W^{\ell}\) are central, so
cancelling their powers in the division ring identifies the
coefficient with the degree-zero center \(C\). Explicitly,
\[
 \begin{aligned}
 C&=\operatorname{span}\{\text{monomials in }D
            \text{ with exponent in }K\},\\
 Z(T)&=C[X^\ell,W^\ell]\subseteq T .
 \end{aligned}
\]
Indeed an element of \(D\) commuting with \(D,X\) commutes with \(W\)
by \(W=X^{-1}U\), and its monomials must be radical. Let \(C_0\)
be the Laurent subalgebra of \(C\) supported on the shared units.
For every central monomial \(M\) of positive zero-frozen degree,
\[
 \chi(M)^\ell=\chi(M^\ell)=0 .
\]
Thus the full-center fiber is imposed by \(c-\chi(c)\) for \(c\in C_0\), these
positive monomials, and \(X^\ell,W^\ell\); there are no additional
central generators. The shared unit torus modulo the \(C_0\)
character is a product \(D_0\) of equal-size matrix algebras by the
finite étale Laurent-center calculation of the preceding chart
proof.

The positive central monomial ideal separates among disjoint
radical groups. In the two groups carrying \(F_1,F_2\), a zero
with all partners units has allowed exponent \(0\). With a zero
partner of exponent \(s_G\), the pair inequality is \(a+s_G<\ell\);
with both triple partners zero it is
\(a+\min(s_{G_1},s_{G_2})<\ell\).
Thus, at fixed spectator exponents \(\mathbf s\), the bounds are
\[
 0\le a<p(\mathbf s),\qquad 0\le b<r(\mathbf s)
\]
for the exponents of \(F_1,F_2\), exactly as in
\eqref{eq:endpoint-bounds}.
Using their unit multiples \(A,B\) preserves these bounds.
The quotient obtained by killing all zero frozens and \(X,W\)
is \(D_0\); in particular the prescribed collision characters
give nonzero fibers, as asserted in (ii).

To compute the remaining ideal, denote the coefficient quotient
by \(\overline D\) after imposing the degree-zero center.
In degree \(|k|\geq\ell\) the entire signed summand is in the
central-power ideal. In degree \(0\leq k<\ell\), put \(m=\ell-k\);
the signed-degree-\(k\) generators are
\[
 \begin{aligned}
 X^\ell W^{\ell-k}
    &=\Bigl(\prod_{j=k}^{\ell-1}U_j\Bigr)X^k,\\
 W^\ell X^{\ell+k}&=N_\ell^- X^k .
 \end{aligned}
\]
The second coefficient is a multiple of the first. In negative
degree \(k=-s\), \(0<s<\ell\), the generators are
\[
 \begin{aligned}
 W^\ell X^{\ell-s}
    &=\Bigl(\prod_{j=s+1}^{\ell}U_{-j}\Bigr)W^s,\\
 X^\ell W^{\ell+s}&=N_\ell^+ W^s .
 \end{aligned}
\]
Again the second coefficient is a multiple of the first. These
are all contributions: the central generators have degrees
\(\ell,-\ell\), so their coefficients in a given degree come
from the corresponding two summands of \(T\). For \(k=0\)
the product is the full-cycle norm, already zero in
\(\overline D\), since \(A^\ell=B^\ell=0\).

Write \(P_m\) for the first displayed product in the appropriate
signed degree. Commutation by shared variables scales all terms
of \(P_m\) by the same scalar, so taking its coefficient ideal
does not mix spectator exponent sectors. The
\((\mathbf s,k)\)-summand has free \(D_0\)-coefficients indexed by
a homogeneous basis of
\[
 H_{\mathbf s,k}=\mathbb C[A,B]/(A^p,B^r,P_m(A,B)),\qquad m=\ell-|k|.
\]
This is a description of coefficient multiplicities, not an
assertion that the \(D_0\)-action is untwisted.

For degree \(j\), dehomogenize at \(A=1\). After scaling and,
if necessary, reversing the order of the roots of \(P_m(1,B)\),
they are \(\gamma_i=\gamma_0\omega^i\) for \(0\leq i<m\).
For a set of exponents \(I\subseteq[0,j]\) the evaluation columns
have entries \(E_{i,b}=\gamma_0^b(\omega^b)^i\), \(b\in I\).
Equal residues give proportional columns; for distinct residues
\(b_1,\ldots,b_d\) with \(d\leq m\), the first \(d\) rows give
\[
 \det (E_{i,b_h})_{\substack{0\leq i<d\\1\leq h\leq d}}
  =\Bigl(\prod_{h=1}^d\gamma_0^{b_h}\Bigr)
      \prod_{1\leq p<q\leq d}(\omega^{b_q}-\omega^{b_p})\ne0 .
\]
Hence \(\operatorname{rank}(E|_I)=\min(m,|I\bmod\ell|)\).
The homogeneous evaluation kernel is \((P_m)_j\), since \(P_m\)
is homogeneous, of degree \(m\), and has \(m\) distinct roots
after dehomogenizing. The degree-\(j\) monomials of the ideal
\((A^p,B^r)\) correspond to \(B_j=[0,j-p]\cup[r,j]\).
As \(m\leq\ell\),
\[
 \begin{aligned}
 \dim(H_{\mathbf s,k})_j
 &=\operatorname{rank}(E|_{[0,j]})-\operatorname{rank}(E|_{B_j})\\
 &=\min(m,j+1)-\min(m,|B_j\bmod\ell|)
   =h_j(\mathbf s,k).
 \end{aligned}
\]

Give \(A,B\) weight two, and \(X,W\) and zero-frozen spectators
weight one. The finite graded fiber has degree-zero quotient
\(D_0\), and its positive part is nilpotent, so that part is its
Jacobson radical \(J\).
The exchange relations and the distinct phases give
\[
 \binom{XW}{WX}
 =
 \begin{pmatrix}1&1\\ \alpha^{-1}&\beta^{-1}\end{pmatrix}
 \binom{A}{B},
 \qquad \beta^{-1}-\alpha^{-1}\ne0 .
\]
Thus \(A,B\in J^2\). All degree-\(d\) monomials are in \(J^d\),
by factoring them into generators of the indicated weights.
Conversely \(J^d\) consists of elements of degree at least \(d\).
Therefore
\[
 J^d=\bigoplus_{b\geq d}(T_\chi)_b .
\]
For a primitive idempotent \(p_\delta\) of \(D_0\), basis
monomials in a fixed sector \((\mathbf s,k,j)\) transport
\(p_\delta\) by the same torus-character permutation: changing
the exponents of \(A,B\) at fixed total \(j\) changes their
commutation by a power of \(AB^{-1}\), which commutes with \(D\).
The sector therefore contributes \(h_j(\mathbf s,k)\) copies of
the same transported simple to the layer
\(|\mathbf s|+|k|+2j\) of \(T_\chi p_\delta\).
Summing over sectors gives \eqref{eq:endpoint-loewy}.
\end{proof}

\begin{remark}
In the unit cancellation of part~(i), the isomorphism to a chart
fiber is already supplied by Theorem~\ref{thm:unit-exchange-patch};
its coefficient corner need not be semisimple.
\end{remark}

\begin{example}\label{ex:node}
Take the endpoint collision with \(F_1,F_2\) the only zero frozen
coordinates. Each has a unit partner, so \(p=r=1\) in
\eqref{eq:endpoint-bounds}. There are no spectator degrees, and
the coefficient quotient for every \(|k|<\ell\) has only its constant
term. Thus
\[
  L_P(z)=1+2\sum_{k=1}^{\ell-1}z^k .
\]
On a seed chart with the same frozen zero set, each zero has
a unit partner and \(L_P(z)=1\). For instance, at \(\ell=3\)
the endpoint projectives have layer lengths \(1,2,2\).
The contrast holds in every size \(n\ge5\).
\end{example}

\begin{proposition}\label{prop:endpoint-center-germ}
At the collision of Example~\ref{ex:node}, the completed local ring
of the full center is
\[
 \mathbb C[[t_1,\ldots,t_{N-3},\mathsf{x},\mathsf{w},u,v]]/(\mathsf{x}\mathsf{w}-u^\ell-v^\ell),
 \qquad N=\binom n2 .
\]
The transverse threefold is an isolated Fermat-type
\(cA_{\ell-1}\) singularity. The exact full-center fiber is
\(M_s(B_\ell)\), where \(s=\operatorname{PIdeg}(T)/\ell\) and
\(B_\ell\) is the algebra of the doubled cycle on \(\mathbb Z/\ell\), with all immediate reversals
and all consistently oriented paths of length \(\ell\) set to zero.
With labels chosen so that \(X\) moves forward, the radical-layer
polynomials are \(C_{ii}(z)=1\) and
\(C_{i+k,i}(z)=z^k+z^{\ell-k}\) for \(1\le k<\ell\).
\end{proposition}

\begin{proof}
Let \(I\) index the shared unit variables and let \(K_{\mathrm{int}}\)
be the integral lift of the initial commutation radical modulo
\(\ell\). Put \(K_I^{\mathrm{int}}=K_{\mathrm{int}}\cap\mathbb Z^I\).
The two frozen radical supports furnish lifts \(k_1,k_2\) of the
standard basis under projection to the \(F_1,F_2\) coordinates.
For \(P_I=\mathbb Z^I\oplus\mathbb N e_{F_1}\oplus\mathbb N e_{F_2}\),
this gives
\[
 \begin{split}
 K_{\mathrm{int}}\cap P_I
     &=K_I^{\mathrm{int}}\oplus\mathbb N k_1\oplus\mathbb N k_2,\\
 C&=C_0[u_0,v_0],\qquad C_0=\mathbb C[K_I^{\mathrm{int}}],\\
 f_1&=g_1u_0^\ell,\qquad f_2=g_2v_0^\ell,\qquad g_i\in C_0^\times,
 \end{split}
\]
where \(u_0,v_0\) are the based monomials of \(k_1,k_2\).
The shared unit lattice has rank \(N-3\), and
\(K_I^{\mathrm{int}}\) has finite index in it. The signed normal
form gives
\[
 Z(T)=C[\mathsf{x},\mathsf{w}]/(\mathsf{x}\mathsf{w}-af_1-yf_2).
\]
At the specified point \(ag_1,yg_2\) are units. Their formal
\(\ell\)-th roots give the stated completed ring by the change
\[
 u=(ag_1)^{1/\ell}u_0,\qquad v=(yg_2)^{1/\ell}v_0.
\]
The transverse equation \(h=\mathsf{x}\mathsf{w}-u^\ell-v^\ell\)
has Jacobian ideal \((\mathsf{x},\mathsf{w},u^{\ell-1},v^{\ell-1})\),
so its singularity is isolated.
For a general section \(v=cu+\alpha\mathsf{x}+\beta\mathsf{w}\),
the Hessian in \(\mathsf{x},\mathsf{w}\) is nondegenerate. The
critical equations in those variables give
\(\mathsf{x},\mathsf{w}\in(u^{\ell-1})\), so the remaining term is
\(-(1+c^\ell)u^\ell+O(u^{2\ell-2})\).
For \(1+c^\ell\ne0\) the splitting lemma gives type \(A_{\ell-1}\).

For the fiber put \(R=\mathbb Z/\ell\),
\(K_I=\ker\Lambda_{II}\) over \(R\), and \(K_0=K\cap R^I\).
For \(k\in K_I\) the identities
\(\langle k,k_1\rangle=\langle k,k_2\rangle=0\) imply
\(\langle k,e_{F_1}\rangle=\langle k,e_{F_2}\rangle=0\),
since the remaining parts of \(k_1,k_2\) are supported on \(I\).
The exponent of \(AB^{-1}\) commutes with the shared variables
and pairs with \(X\) by a unit modulo \(\ell\). Subtracting its
two frozen coordinates times \(k_1,k_2\) preserves its pairings
and yields an exponent supported on \(I\). Consequently
\[
 0\longrightarrow K_0\longrightarrow K_I
   \xrightarrow{\ \langle\,\cdot,e_X\rangle\ }R\longrightarrow0,
 \qquad
 D_0\simeq\bigoplus_{i\in R}M_s(\mathbb C),\qquad
 \sigma(e_i)=e_{i+1}.
\]
Here \(D_0\) denotes the shared unit algebra after specialization,
and \(e_i\) are its central block idempotents; the last equality
uses the surjectivity in the exact sequence. Since the central
lifts kill \(F_1,F_2\), both exchange monomials vanish in the
fiber. Write \(A=u_0a_0\) before specialization. Conjugation by
\(W\) on shared units is \(\operatorname{Ad}(A)\circ\sigma^{-1}\).
Thus \(\widetilde W=a_0^{-1}W\) induces \(\sigma^{-1}\), and
\[
 Xe_i=e_{i+1}X,\qquad
 \widetilde We_i=e_{i-1}\widetilde W,\qquad
 X\widetilde W=\widetilde WX=X^\ell=\widetilde W^\ell=0 .
\]
Choose the same diagonal matrix unit \(p_i\) in each block by
transport under \(\sigma\), and set \(e=\sum_i p_i\).
The signed normal form with \(p=r=1\) in
\eqref{eq:endpoint-bounds} gives the corner basis consisting of
\(p_i,\ X^kp_i,\ \widetilde W^kp_i\) for \(i\in R\) and
\(1\le k<\ell\) (with \(p_i\) included once).
The path relations above give the same basis for \(B_\ell\),
so \(eT_\chi e\cong B_\ell\). Summing corresponding matrix
units across the blocks gives commuting matrix units for
\(T_\chi\cong M_s(B_\ell)\), with size given by
\[
 s^2=\frac{\ell^{N-3}}{|K_I|}
      =\frac{\ell^{N-3}}{\ell|K_0|}
      =\frac{\ell^{N-2}}{|K|}
      =\frac{\operatorname{PIdeg}(T)^2}{\ell^2}.
\]
Here compatibility forces the \(X\)-coordinate of \(K\) to
vanish over \(R\), and the lifts \(k_1,k_2\) give
\(|K|=\ell^2|K_0|\).
Finally the radical is the arrow ideal: its quotient is
\(\mathbb C^\ell\) and it is nilpotent. From \(i\) to \(i+k\)
there are exactly two nonzero monotone paths, of lengths \(k\)
and \(\ell-k\), for \(0<k<\ell\); for \(k=0\) only the idempotent
remains. This gives the stated polynomials.
\end{proof}

\begin{remark}
The fixed-frozen nodal slice sets \(u=v=0\); it should not be
confused with the transverse germ of the full center.
For the symplectic cohomology of Milnor fibers of the latter
Fermat family, see \citet[Theorem C]{APZ}.
\end{remark}

\begin{proposition}\label{prop:endpoint-thickening}
Under the assumptions of Example~\ref{ex:node}, fix the smooth
parameters of the full center at the collision. The resulting
completed algebra is a matrix algebra over the completed doubled
\(\ell\)-cycle \(\mathcal B_\ell\) over \(\mathbb C[[u,v]]\).
If \(a_i:i\to i+1\) and \(b_i:i+1\to i\) are its arrows
(indices modulo \(\ell\)), its relations are
\[
 b_i a_i=(u+\omega^i v)e_i,\qquad
 a_i b_i=(u+\omega^i v)e_{i+1},
\]
for a primitive \(\ell\)-th root \(\omega\) chosen to match the
cyclic labeling. The sums of the two
kinds of oriented full cycles, denoted by \(\mathsf{x},\mathsf{y}\), are central
and satisfy \(\mathsf{x}\mathsf{y}=u^\ell+v^\ell\). Completion is at the ideal
\((\mathsf{x},\mathsf{y},u,v)\) in the center.

At skew root order \(\ell=3\), \(\mathcal B_3/(\mathsf{x},\mathsf{y},uv)\) is
isomorphic to the three-vertex basic algebra of the generalized
reduction of the subregular \(A_2\) Frobenius block at a regular
Harish--Chandra character of Theorem~\ref{thm:a2-cartan}.
Further quotienting by \((u,v)\) gives the basic algebra of the
exact Harish--Chandra quotient. The quantum-group root order is
independent of the skew root order and retains the restrictions
of that theorem.
\end{proposition}

\begin{proof}
Use the central lifts \(u_0,v_0\) from the preceding proof and write
\(D=Q[u_0,v_0]\), with \(Q\) the shared unit quantum torus. Put
\[
 A=u_0 a_0,\qquad B=v_0 b_0,\qquad
 S=b_0a_0^{-1},\qquad Y=Wa_0^{-1}.
\]
The units \(a_0,b_0\) have the same conjugation on \(Q\), and
\(a_0\) is a \(\sigma\)-eigenvector. Hence \(S\in Z(Q)\) and
\[
 Yd=\sigma^{-1}(d)Y,qquad
 XY=u_0+v_0S,\qquad YX=\sigma^{-1}(u_0+v_0S)
 \quad(d\in D).
\]
Also \(Y^\ell=W^\ell a_0^{-\ell}\), since the scalar
from commuting the factors is an \(\ell\)-th root of unity
to the power \(\ell(\ell-1)/2\). The lattice calculation above
gives \(Z(Q)=C_0[S^{\pm1}]\), with \(S^\ell\in C_0^\times\)
and \(\sigma(S)/S\) a primitive root. Over the complete local
ring of \(C_0\), choose \(\rho\) with \(\rho^\ell=S^\ell\)
and put \(u=u_0\), \(v=\rho v_0\). The center of \(Q\) splits
according to the distinct roots of \(S^\ell-\rho^\ell\).
Over the completed coefficient ring, each central factor of
\(Q\) is split Azumaya, since the residue field is \(\mathbb C\).
Write \(e_i\) for the factors ordered by \(\sigma(e_i)=e_{i+1}\).
For matrix units \(E^{(0)}_{jk}\) in the first factor, put
\[
 E_{jk}=\sum_{i=0}^{\ell-1}\sigma^i(E^{(0)}_{jk}).
\]
Then \(\sigma(E_{jk})=E_{jk}\), so \(X,Y\) commute with
these matrix units. The corner cut out by \(E_{11}\) has
vertices \(e_iE_{11}\); the two shifts are induced by \(X,Y\).
Choosing \(\omega\) according to this order, their
backtracking coefficient at bond \(i\) is \(u+\omega^iv\).

In the path algebra with these reversal relations, a path
reduces to a scalar multiple of a monotone path. For each
starting vertex there is one monotone path in each signed
degree, with the degree-zero path the vertex idempotent.
Their linear independence before completion follows from
the signed normal form of Theorem~\ref{thm:endpoint}.
Taking the central completion and then fixing the smooth
parameters therefore gives the stated presentation.
In particular
\[
 \mathsf{x}=X^\ell,\qquad \mathsf{y}=Y^\ell,\qquad
 \mathsf{x}\mathsf{y}=\prod_{i=0}^{\ell-1}(u+\omega^i v)=u^\ell+v^\ell.
\]

Now let \(\ell=3\). Set \(\xi_i=\omega^i u+\omega^{2i}v\),
\(0\le i<3\). This linear map identifies \(\mathbb C^2\) with the trace-zero plane, and its elementary symmetric
functions are
\[
 \sum_i\xi_i=0,\qquad
 \sum_{i<j}\xi_i\xi_j=-3uv,\qquad
 \xi_0\xi_1\xi_2=u^3+v^3 .
\]
Thus in the quotient by \((\mathsf{x},\mathsf{y},uv)\) the scalar ring
\(\mathbb C[[u,v]]/(uv,u^3+v^3)\) is the Harish--Chandra coefficient
algebra in Theorem~\ref{thm:a2-cartan}. Replacing \(b_i\) by
\(\omega^i b_i\) makes both reversals at bond \(i\) equal to
\(\xi_i\), without changing the zero cycle relations.
A length-one arc is annihilated by the complementary product
of length two, and a length-two arc by the complementary bond.
The signed normal form modulo the central cycle ideal shows
that these are precisely the coefficient ideals: in degree
\(0<|k|<3\), the only contributing cycle multiple uses the
opposite complementary degree \(3-|k|\), and in degree zero
the product of all three bonds already vanishes. The resulting
three-vertex presentation is the contracted presentation in
the proof of Theorem~\ref{thm:a2-cartan}. Its dimension is
\(3\cdot6+6\cdot(4+2)=54\).
Finally the ideals \((u,v)\) and \((\xi_0,\xi_1,\xi_2)\)
correspond under the same linear change of coordinates.
\end{proof}

\begin{remark}
The algebra \(\mathcal B_\ell\) is the known noncommutative
crepant resolution of the threefold with equation
\(\mathsf{x}\mathsf{y}=\prod_i(u+\omega^i v)\), in the doubled-cycle presentation
of \citet[Section~2.2 and Corollaries~4.6--4.7]{EL} and the
smooth-factor case of \citet[Corollary~5.4]{IW}.
Here the new identification is its occurrence from the
quantum endpoint exchange and the explicit central
quotient giving the generalized block.
\end{remark}

\subsection{Projective layers at all single-coordinate nodes}
\label{sec:general-nodal-layers}

In the three terminal types of Proposition~\ref{prop:single-mutable-coverage}, the two exchange monomials may involve more than
two vanishing frozen variables. The formula
\eqref{eq:endpoint-rank} extends to this situation, but the
forbidden exponents need not form the two intervals occurring
there. We describe them directly in terms of the radical lattice.

\begin{theorem}\label{thm:all-single-nodal-layers}
Let \(n\geq6\) and let \(\ell>1\) be odd. Invert the powers of all
initial mutable variables except \(X\), let \(W\) be its exchange
partner, and write \(XW=A+B\) in the based-torus normalization.
Let \(\chi\) be a full-center character for which
\[
 \chi(X^\ell)=\chi(W^\ell)=\chi(A^\ell)=\chi(B^\ell)=0.
\]
Let \(S\) be the set of zero frozen variables, and write
\(S_A,S_B\subseteq S\) for the supports of the zero frozen
factors in \(A,B\). These sets are disjoint and nonempty.
In these three terminal types,
each factor occurs to the first power, and
\(\{|S_A|,|S_B|\}\) is either \(\{1,1\}\) or \(\{1,2\}\).

Write \(\pi_S\) for exponent projection onto \(S\) and
\(\bar u\) for reduction modulo \(\ell\). Call \(e\in\mathbb N^S\) standard if
\[
 e\notin\bigcup_{\substack{u\in\mathbb Z^N,\ \bar u\in K,\ u_X=0\\
                         0\ne\pi_S(u)\in\mathbb N^S}}
                \bigl(\pi_S(u)+\mathbb N^S\bigr).
\]
Give each coordinate in a singleton \(S_A\) or \(S_B\) weight
two, and every other coordinate in \(S\) weight one; write
\(\mathrm{wt}(e)\) for the resulting weight. Sum over
\[
 \mathbf s\in\{0,\ldots,\ell-1\}^S,\qquad
 \min_{i\in S_A}s_i=\min_{i\in S_B}s_i=0 .
\]
For \(|k|<\ell\) and \(0\leq j\leq 2\ell-2\), let
\[
 \begin{split}
 R_j(\mathbf s)&=
 \{\,b\bmod\ell:0\leq b\leq j,\
   \mathbf s+(j-b)\mathbf 1_{S_A}+b\mathbf 1_{S_B}
                   \text{ is not standard}\,\},\\
 h_j(\mathbf s,k)&=\min(\ell-|k|,j+1)-
                       \min(\ell-|k|,|R_j(\mathbf s)|).
 \end{split}
\]
Then the unlabelled Loewy polynomial of every projective cover is
\begin{equation}\label{eq:all-single-nodal-layers}
 \sum_{\mathbf s}\sum_{|k|<\ell}\sum_{j=0}^{2\ell-2}
       h_j(\mathbf s,k)z^{\,\mathrm{wt}(\mathbf s)+2j+|k|}.
\end{equation}
In each summand, the simple label is obtained by the common
torus-character permutation of the monomials in that summand.
The fiber is not Azumaya.
\end{theorem}

\begin{proof}
Invert also the powers of the frozen variables outside \(S\).
Let \(D\) be the mixed coefficient torus, with exactly the
coordinates in \(S\) noninverted. By compatibility, the two
exchange monomials commute, \(AB^{-1}\) commutes with \(D\)
in the quotient division ring, and, for
\(\sigma=\operatorname{Ad}_X\),
\[
 \sigma(A)=\alpha A,\qquad \sigma(B)=\beta B,\qquad
 \operatorname{ord}(\alpha/\beta)=\ell .
\]
These are precisely the identities used in the signed-degree
comparison in Theorem~\ref{thm:endpoint}; the cardinalities
of \(S_A,S_B\) play no role in that comparison. If \(T\)
denotes the resulting partial localization, it gives
\[
 \begin{aligned}
 T_k&=\begin{cases} DX^k,& k\geq0,\\ DW^{-k},&k<0,\end{cases}
 & Z(T)&=C[X^\ell,W^\ell],\\
 C&=\operatorname{span}\{M(u)\in D:u\in\mathbb Z^N,\ \bar u\in K,\ u_X=0\}.
 \end{aligned}
\]
Here \(M(u)\) denotes a based monomial. The character specializes
the degree-zero Laurent part of \(C\) to scalars and kills every
central monomial with positive \(S\)-degree: its \(\ell\)-th
power has zero value. Let \(D_0\) be the shared unit torus
modulo that Laurent character. As in Theorem~\ref{thm:endpoint},
\(D_0\) is a product of matrix algebras. After imposing the
character on \(C\) the coefficient module has one free left
\(D_0\)-summand per standard exponent. The divisibility test
above follows by multiplying each positive central monomial
by arbitrary unit exponents and nonnegative \(S\)-exponents.

Write
\[
 e=\mathbf s+a\mathbf1_{S_A}+b\mathbf1_{S_B},\qquad
 a=\min_{i\in S_A}e_i,\quad b=\min_{i\in S_B}e_i .
\]
For each \(\mathbf s\) choose the ordered monomial \(F^{\mathbf s}\)
in the zero frozen variables. In degree \(k\), the indicated
basis consists, up to left unit factors, of terms
\(F^{\mathbf s}A^aB^bX^k\) for \(k\geq0\), and
\(F^{\mathbf s}A^aB^bW^{-k}\) for \(k<0\).
Multiplication by \(A\) or \(B\) preserves \(\mathbf s\).
As conjugation by coefficient monomials scales \(A\) and \(B\)
by the same factor, the central-power relations are computed
separately in each \((\mathbf s,k)\)-sector.

In particular \(|k|\geq\ell\) gives zero. For \(|k|<\ell\)
put \(m=\ell-|k|\). The calculation with the two adjacent
Laurent charts gives in signed degree \(k\) one coefficient
relation, the product \(P_m(A,B)\) of \(m\) consecutive
shifted exchange binomials; when \(k=0\) it is the full-cycle
norm \(A^\ell+B^\ell\). At each \(\mathbf s\), the resulting
scalar quotient is by this relation and by the monomials
\(A^aB^b\) with nonstandard exponent \(e\). To see the
dimension in degree \(j=a+b\), set \(A=1\) and evaluate in
the \(m\) roots of \(P_m(1,B)\),
\(\gamma_i=\gamma_0(\alpha/\beta)^{\pm i}\), \(0\leq i<m\).
The kernel is the degree-\(j\) part of \((P_m)\). Writing
\(\omega=(\alpha/\beta)^{\pm1}\), the column of \(A^{j-b}B^b\)
is \(\bigl(\gamma_0^b(\omega^b)^i\bigr)_{i=0}^{m-1}\).
For a set \(I\) of such exponents the Vandermonde determinant gives
\[
 \operatorname{rank}
    \bigl(\gamma_0^b(\omega^b)^i\bigr)_{0\leq i<m,\,b\in I}
     =\min(m,|I\bmod\ell|).
\]
Subtract the rank contributed by the forbidden columns from the
rank of all \(j+1\) columns to obtain \(h_j(\mathbf s,k)\).
For \(j>2\ell-2\), every exponent has \(a\geq\ell\) or
\(b\geq\ell\), so every column is forbidden. The calculation
uses the order \(\ell\) of \(\omega\), not primality of \(\ell\).

Give \(X,W\) degree one, use the indicated weights on \(S\)
and degree zero on the shared units. Then \(A,B\) each have
degree two. The resulting grading on the finite fiber
\(T_\chi\) has semisimple degree-zero part \(D_0\); hence its
positive part is the Jacobson radical \(J\).
Moreover
\[
 \binom{XW}{WX}=
 \begin{pmatrix}1&1\\ \alpha^{-1}&\beta^{-1}\end{pmatrix}
 \binom{A}{B}
 \quad\Longrightarrow\quad A,B\in J^2 .
\]
Every weight-two frozen generator is a unit multiple of
\(A\) or \(B\); thus all algebra generators of positive
weight \(d\) lie in \(J^d\). This yields
\[
 J^d=\bigoplus_{i\geq d}(T_\chi)_i .
\]
At fixed \((\mathbf s,k,j)\) changing \(a,b\) with \(a+b=j\)
does not change the induced conjugation on \(D_0\), since
\(AB^{-1}\) commutes with \(D\). The operators \(F^{\mathbf s}\),
\(A,B,X,W\) used for this conjugation are taken in the quotient
division ring before specialization; the induced automorphisms
preserve the fixed character of the Laurent part of \(C\).
Consequently the sector is \(h_j(\mathbf s,k)\) copies of the
regular \(D_0\)-bimodule with the same twist on the right.
Multiplying by a primitive idempotent of \(D_0\) gives the
claimed composition length and transported simple label in
weight \(\mathrm{wt}(\mathbf s)+2j+|k|\).
Finally \(h_0(\mathbf0,1)=1\), so \(J/J^2\ne0\).
\end{proof}

\subsection{Simultaneous unit exchanges and the deep locus}
\label{subsec:odd-pair-open}

We now allow several mutable coordinates to vanish at once. In contrast
to a node with two vanishing exchange monomials, a point may lie outside
every seed chart even though all the exchange monomials under consideration
are units. The following family exhibits this distinction in every even
size \(n=2m\ge6\).

Let \(\mathcal I=\{(i,j):1\le i<j\le 2m-1,\ i,j\text{ odd}\}\) and
\(r=\binom m2\). Write \(X_a\) for the initial variable at \(a\in\mathcal I\)
and \(W_a\) for its mutation. Put \(\mathsf x_a=X_a^\ell\) and
\(\mathsf w_a=W_a^\ell\). Let \(U_0\subset\operatorname{Spec}Z_0\)
be the open set on which all initial powers outside \(\mathcal I\) are
units, and \(U\subset\operatorname{Spec}Z(A)\) its inverse image.
On \(U\), denote by \(\Delta\) the reduced locus outside every seed chart.

\begin{theorem}\label{thm:odd-pair-open}
The algebra is Azumaya throughout \(U\). A point of \(U\) is outside
every seed chart if and only if
\[
 \mathsf x_a=\mathsf w_a=0\qquad(a\in\mathcal I).
\]
Every other point is covered by a seed obtained using at most \(2r\)
mutations. The locus \(\Delta\) and its image in \(U_0\), with reduced
structure, are tori of dimension \(m^2\) over \(\mathbb C\); their
map is finite \'etale of degree \(\ell^m\).
\end{theorem}

\begin{proof}
For a position \((i,j)\), use \(e_{ij}\) for its standard exponent.
At \((i,j)\in\mathcal I\), the columns from
\eqref{eq:all-exchange-columns} are
\begin{equation}\label{eq:odd-pair-columns}
\begin{split}
 b_{ij}={}&e_{i,j+1}-e_{i,j-1}+e_{i+1,j}-e_{i+1,j+1}\\
 &+\mathbf1_{\{i>1\}}(e_{i-1,j-1}-e_{i-1,j})
 +\mathbf1_{\{j=i+2\}}
    (e_{i,i+1}-\mathbf1_{\{i>1\}}e_{i-1,i}).
\end{split}
\end{equation}
In particular the selected positions are pairwise nonadjacent,
and their exchange monomials use only variables outside \(\mathcal I\).
Join consecutive selected positions in each row and column. This
graph is connected: each row leads to the last odd column, which
is itself connected. For a horizontal edge from \((i,j)\) to
\((i,j+2)\) the intervening row of the exchange matrix is \((i,j+1)\);
for a vertical edge to \((i+2,j)\) it is \((i+1,j)\).
By \eqref{eq:odd-pair-columns}, each such row has entries \(1,-1\) on
the two endpoints (up to a common sign) and zero on the other
selected columns. These intervening positions are mutable and distinct.
Summing \eqref{eq:odd-pair-columns}, the horizontal differences telescope
in each odd row. In even row \(i\) the two adjacent odd rows contribute
\(e_{i,i+1}-e_{i,n}\); the near-diagonal terms cancel after adding the
last correction in \eqref{eq:odd-pair-columns}. Thus
\begin{equation}\label{eq:odd-pair-sum}
 \sum_{a\in\mathcal I}b_a
 =\sum_{i=1}^{n-2}(-1)^{i+1}e_{i,n}+e_{n-2,n-1}.
\end{equation}
This is supported on frozen rows. Take the edge rows of a spanning tree and adjoin row \((1,n)\).
In the resulting \(r\times r\) minor replace one column by the
column sum. Up to row and column permutations its determinant is
\[
 \det\begin{pmatrix}E&0\\ *&1\end{pmatrix}
 =\det E=\pm1,
\]
where \(E\) is the signed tree incidence matrix with the replaced
column deleted. Thus the selected columns span a primitive
sublattice of the shared exponent lattice.

Consider the Frobenius exchanges \(\mathsf x_a\mathsf w_a=
A_a^\ell+B_a^\ell\), where \(X_aW_a=A_a+B_a\), with the summands
ordered to give the column \(b_a\) as their exponent difference.
Let \(d\) have weight one on \(\mathcal I\) and zero elsewhere.
For each mutable index \(k\), skew symmetry and
\eqref{eq:odd-pair-sum} give
\[
 (d^{\mathsf T}B)_k=\sum_{a\in\mathcal I}b_{ak}
       =-\Bigl(\sum_{a\in\mathcal I} b_a\Bigr)_k=0.
\]
After mutation at \(k\), the weights are
\(d'_j=d_j\) for \(j\ne k\) and
\(d'_k=-d_k+\sum_j[b_{jk}]_+d_j\). This invertible change
preserves homogeneity of the exchanges, and the frozen weights
remain zero. In particular every seed contains a mutable power
of nonzero weight. The grading on \(U_0\) sends
\((\mathsf x_a,\mathsf w_a)\) to
\((t\mathsf x_a,t^{-1}\mathsf w_a)\) and fixes all shared powers;
the presentation below shows that these are generators of
\(\mathcal O(U_0)\). Consequently a point with all pairs zero is
fixed. For a nonzero weight \(d'\) and corresponding mutable
power \(y'\) in any seed, its value at such a point satisfies
\((t^{d'}-1)y'=0\) for every \(t\in\mathbb C^\times\).
Thus \(y'=0\), excluding that seed chart.

Conversely, swap every \(\mathsf x_a=0\) for its partner if that
partner is nonzero. These swaps commute. If the remaining set of
collisions is nonempty and proper, choose a graph edge from this set
to its complement. Its intervening variable \(z\) is a unit and is
adjacent to just one unresolved exchange \(xv=za+b\). The exchange
at \(z\) has the form \(zz'=c+xd\), with \(a,b,c,d\) units.
By the two-step calculation in
Section~\ref{subsec:single-mutable-coverage}, mutation at \(z\) and
then at \(x\) gives, at the collision,
\[
 z'=c/z,\qquad u=(b/g)d/z,\qquad g=\gcd(b,c).
\]
Both values are units. If \(J\) is the set of unresolved positions
before the step and \(a\) the endpoint removed, the exchange matrix
satisfies \(b_{zk}=b_{ak}=0\) for \(k\in J\setminus\{a\}\).
The matrix mutation formula leaves each column at such \(k\)
unchanged in both steps; the intervening row cannot be needed
again. If there were \(h<r\) initial collisions, there are at most
\(r-h\) swaps and \(2h\) further mutations, hence at most \(2r\).
The criterion in the full center is the pullback of this criterion,
because seed powers belong to \(Z_0\).

For the algebra calculation let \(M\) be the shared exponent lattice,
\(D\) its quantum torus, and set
\[
 t_a=A_aB_a^{-1},\qquad \rho_a=t_a^\ell,\qquad
 \sigma_a=\operatorname{Ad}_{X_a}|_D,\qquad
 D_*=\bigcap_{a\in\mathcal I}D^{\sigma_a}.
\]
Compatibility gives \(t_a\in Z(D)^\times\) and
\(\sigma_b(t_a)=\omega_a^{\delta_{ab}}t_a\), with \(\omega_a\) of
order \(\ell\). Let \(\mathcal L\) be the localization over \(U_0\).
In the common fraction ring the \(2^r\) seeds obtained by independent
swaps have Laurent degrees in \(\mathbb Z^{\mathcal I}\), assigning
\(\deg X_a=e_a,\ \deg W_a=-e_a,\ \deg D=0\).
For \(\mathbf d\in\mathbb Z^{\mathcal I}\) let \(S_{\mathbf d}\)
be the ordered product with factors \(X_a^{d_a}\) for \(d_a\ge0\)
and \(W_a^{-d_a}\) for \(d_a<0\). In the chart using precisely these
variables, the degree-\(\mathbf d\) subspace is \(DS_{\mathbf d}\).
Let \(\mathcal T_{\varepsilon}\) denote the localized torus for
the independent swaps specified by \(\varepsilon\in\{0,1\}^{\mathcal I}\).
The Laurent inclusions and \(DS_{\mathbf d}\subset\mathcal L\) give
\[
 \bigoplus_{\mathbf d}DS_{\mathbf d}\ \subseteq\ \mathcal L
 \ \subseteq\ \bigcap_{\varepsilon\in\{0,1\}^{\mathcal I}}\mathcal T_\varepsilon
 \ \subseteq\ \bigoplus_{\mathbf d}DS_{\mathbf d}.
\]
All tori in this intersection have \(D\) as the shared torus. The same argument
in the Frobenius algebra gives the commutative presentation
\[
 \mathcal O(U_0)
 \simeq
 \frac{\mathbb C[\ell M][\mathsf x_a,\mathsf w_a:a\in\mathcal I]}
 {(\mathsf x_a\mathsf w_a-A_a^\ell-B_a^\ell:a\in\mathcal I)} .
\]
This also shows that the grading used above acts on \(U_0\).

Order the selected indices and write \(X_aX_b=\zeta^{c_{ab}}X_bX_a\).
For \(a<b\), choose \(u_{ab}\) with \(\omega_a^{u_{ab}}\zeta^{c_{ab}}=1\).
Then
\[
 \widehat X_b=\Bigl(\prod_{a<b}t_a^{u_{ab}}\Bigr)X_b,\qquad
 \widehat X_a\widehat X_b=\widehat X_b\widehat X_a .
\]
The conjugations on \(D\) are unchanged. Set
\(\widehat Z_a=\widehat X_a^{-1}(1+t_a)\). These replace \(X_a,W_a\)
up to units of \(D\). Set \(C=Z(D_*)\), and let \(Q_a\) be the
rank-one algebra over \(C\) with variables \(t_a,\widehat X_a,\widehat
Z_a\) and relations as in Theorem~\ref{thm:unit-exchange-patch}.
The pairs commute across distinct indices. The simultaneous
eigenspace decomposition and multiplication give
\begin{equation}\label{eq:odd-pair-factorization}
 D=\bigoplus_{0\le d_a<\ell}D_*\prod_a t_a^{d_a},
 \qquad
 \mathcal L\simeq D_*\otimes_C\bigotimes_{a\in\mathcal I}{}_C Q_a .
\end{equation}
Indeed, the ordered signed monomials in the rank-one factors, with
the indicated powers of the \(t_a\), recover the signed basis over
\(D\) after the unit changes. Each \(Q_a\) is Azumaya of degree
\(\ell\) over its center by Theorem~\ref{thm:unit-exchange-patch}.
The invariant algebra \(D_*\) is a quantum torus, hence Azumaya
over \(C\). Write \(S_a=Z(Q_a)\). Over an \'etale cover splitting these
factors and \(D_*\), the right-hand side of
\eqref{eq:odd-pair-factorization} becomes a matrix algebra over
\[
 S=\bigotimes_{a\in\mathcal I} {}_C S_a,\qquad
 S_a=C[\widehat X_a^\ell,\widehat Z_a^\ell]/
       (\widehat X_a^\ell\widehat Z_a^\ell-1-\rho_a).
\]
Faithfully flat descent gives \(Z(\mathcal L)=S\) and proves
the Azumaya assertion.

Let \(M_*\subset M\) be the exponent lattice of \(D_*\), and let
\(\Gamma\) be the full seed radical lattice restricted to \(M\).
Since \(D\) is generated by \(D_*\) and the central units \(t_a\),
a monomial of \(D_*\) is central in \(D_*\) exactly when it belongs
to the full seed radical. Thus \(C=\mathbb C[\Gamma]\) in based
monomials. Compatibility forces radical exponents to be divisible
by \(\ell\) at mutable positions. The radical computation of
Section~\ref{subsec:radical} therefore gives
\(\Gamma/\ell M\simeq(\mathbb Z/\ell)^m\). Write
\(L_\rho=\sum_a\mathbb Z\ell b_a\).
If an element of \(\Gamma\) is in \(\mathbb Q L_\rho\), by the primitive
minor it is \(\sum_a k_ab_a\) with integral \(k_a\). Pairing with the
selected exponents gives \(2k_a\equiv0\bmod\ell\), so \(L_\rho\) is primitive
in \(\Gamma\). It is primitive in \(\ell M\) by the same minor.
On the simultaneous collision locus all rank-one powers are zero
and \(\rho_a=-1\). After translation by a point of the
Laurent-character equations, the map of reduced loci is dual to
\[
 0\longrightarrow \ell M/L_\rho
 \longrightarrow \Gamma/L_\rho
 \longrightarrow(\mathbb Z/\ell)^m\longrightarrow0 .
\]
Both quotient lattices are free of rank
\(\operatorname{rk}M-r=\binom n2-2r=m^2\).
The resulting isogeny of tori is finite \'etale of degree \(\ell^m\).
\end{proof}

\begin{corollary}\label{cor:odd-pair-brauer}
The Brauer index of \(A\) at the generic point of \(U\) is \(\ell^{2r}\).
At the generic point of \(\Delta\) it is \(\ell^r\).
Every closed full-center fiber on \(U\) is \(M_{\ell^{2r}}(\mathbb C)\).
\end{corollary}

\begin{proof}
The generic fiber of the initial quantum torus over its center is a
division algebra: it is a finite-dimensional central localization
of a domain. Its degree is \(\ell^{2r}\) by \eqref{eq:pidegree}.
On \(\Delta\) each rank-one factor in
\eqref{eq:odd-pair-factorization} specializes to the constant
matrix algebra \(M_\ell(\mathbb C)\), since its powers are zero
and its ratio power is \(-1\). Thus
\[
 \mathcal L|_\Delta\simeq M_{\ell^r}\bigl(D_*/(\rho_a+1:a\in\mathcal I)\bigr).
\]
The quotient in parentheses is a quantum torus on \(M_*/L_\rho\):
\(L_\rho\) is primitive in \(M_*\) by the same pairing argument as for
\(\Gamma\), and a lattice complement realizes the quotient by the
specified central characters. It is therefore a domain. The index
calculation is
\[
 [M:M_*]=\ell^r,\qquad [M:\Gamma]=\ell^{\binom n2-r-m}
 =\ell^{3r},\qquad [M_*/L_\rho:\Gamma/L_\rho]=\ell^{2r}.
\]
Here the first equality follows from the independent characters
of the \(t_a\) under the \(\sigma_b\). Consequently the generic
central localization of the residual torus is a division algebra
of degree \(\ell^r\). The assertion at closed points follows from
Theorem~\ref{thm:odd-pair-open} and the degree of the ambient algebra.
\end{proof}

\section{Block-diagonal fibers outside the seed charts}\label{sec:block}

Write \(n=2m\) and \(I_i=\{2i-1,2i\}\). For a character
\(\chi_0:Z_0\to\mathbb C\) consider the following power fiber:
\[
 \begin{gathered}
 B_0=A/\ker(\chi_0)A,\qquad d_i=x_{2i-1,2i},\\
 \chi_0(d_i^\ell)=c_i\ne0,\qquad
 \chi_0(x_{ab}^\ell)=0\quad(a\in I_i,\ b\in I_j,\ i<j).
 \end{gathered}
\]
Choose \(\rho_i^\ell=c_i\), put \(G=(\mathbb Z/\ell)^m\), and
write \(u^\gamma\) for the basis element of the group algebra
\(\mathbb Z[G]\) attached to \(\gamma\in G\).
Let \(\varepsilon_i\) denote the coordinate vectors of \(G\).

\begin{theorem}\label{thm:block}
The simple \(B_0\)-modules are one-dimensional and indexed by
\(\alpha\in G\), with \(d_i\) acting by \(\rho_i\zeta^{\alpha_i}\).
Their projective multiplicities are
\begin{equation}\label{eq:block-product}
 \sum_{\beta\in G} C_{\beta\alpha}(t)u^{\beta-\alpha}
 =
 \prod_{i<j}
    \left(\sum_{h=0}^{\ell-1}
       t^{(j-i)h}u^{h(\varepsilon_i-\varepsilon_j)}\right)^4.
\end{equation}
The \(\ell\) full-center points above \(\chi_0\) separate
the simple labels by \(\sum_i\alpha_i\). Each corresponding
fiber has \(\ell^{m-1}\) one-dimensional simples and lies
outside every seed chart. Formula \eqref{eq:block-product}
computes the layers over \(Z_0\), prior to the remaining
full-center relations.
\end{theorem}

\begin{proof}
For \(a<b\) put \(\delta_k(x_{ab})=1\) if \(a\le2k<b\),
and zero otherwise. Extend additively to a PBW monomial
\(x^v\) and set \(w(x^v)=\sum_k\delta_k(x^v)\).
For \(r<s<t<u\) the three pairings satisfy, componentwise,
\[
 \delta(x_{rs}x_{tu})\le
 \delta(x_{rt}x_{su})=\delta(x_{ru}x_{st}).
\]
In the fiber, each nonzero power reduction occurs on a \(d_i\),
which has weight zero. Hence for the subspaces
\(V_b=\langle x^v:w(x^v)\ge b,\ 0\le v_i<\ell\rangle\),
\[
 \begin{gathered}
 V_bV_c\subseteq V_{b+c},\qquad
 B_0/V_1=\bigotimes_i\mathbb C[d_i]/(d_i^\ell-c_i),\\
 V_b=0\quad\bigl(b>4(\ell-1)\textstyle\sum_{i<j}(j-i)\bigr).
 \end{gathered}
\]
The quotient is semisimple, so \(J(B_0)=V_1\).
If \(a\in I_i,b\in I_j\), \(j>i+1\), take
\(r<s\) in \(I_{i+1}\). The crossing relation gives
\[
 x_{ab}=(\zeta-\zeta^{-1})^{-1}
              [x_{as},x_{rb}]\,d_{i+1}^{-1}\in J(B_0)^{j-i}
\]
by induction on \(j-i\). Thus \(V_b=J(B_0)^b\) in all degrees.

In the associated graded algebra, the commutation with \(d_k\)
of the symbol of a crossing entry \(x_{ab}\) is
\[
 d_k\overline{x}_{ab}
  =\zeta^{\delta_{ki}-\delta_{kj}}\overline{x}_{ab}d_k
       \quad(a\in I_i,b\in I_j).
\]
For \(k=i,j\) the entries share a vertex; for \(i<k<j\)
they are nested; for exterior \(k\) the two correction terms
have larger weight. These are all the cases. Consequently the
graded component of a crossing exponent vector \(h=(h_{ab})\)
is a free rank-one right module over the degree-zero algebra
with a generator of weight and label shift
\[
 \left(\sum_{i<j}\sum_{a\in I_i,\, b\in I_j}(j-i)h_{ab},
 \quad
 \sum_{i<j}\sum_{a\in I_i,\, b\in I_j}h_{ab}
                      (\varepsilon_i-\varepsilon_j)\right).
\]
The PBW basis proves freeness: move the degree-zero factors
to the right using the displayed scalar relation in the graded
algebra. Multiplication by the idempotent for \(\alpha\) on the
right leaves one copy of the simple with the resulting shifted
label. Summing over \(0\le h_{ab}<\ell\) gives
\eqref{eq:block-product}.

Let \(\overline H\) denote the image of \(H\) in \(B_0/J(B_0)\).
The row weight of \(H\) is \((1,\ldots,1)\).
Among the weight-zero PBW monomials for the cut filtration,
the unique monomial of this row weight is \(d_1\cdots d_m\),
its leading term, so
\[
 \overline H=c\,d_1\cdots d_m\quad(c\in\mathbb C^\times).
\]
The remaining supported frozen coordinates are \(B\) and
the paired \(F_r\). For each of these, the
leading chain crosses a cut, and the coordinate
generically quasi-commutes with every \(d_i\) with the phase
of that leading exponent. In degree zero, any nonzero image
would commute with every \(d_i\), contrary to the nonzero
integer shift of a cut-crossing leading word. This works
in the generic cut quotient, a polynomial ring in the
\(d_i\), before imposing any power relations. Hence \(H\) is the only unit on the
radical supports of Theorem~\ref{thm:intro-center}. By
\eqref{eq:center-carry} the residue monomials involving any of those coordinates
are nilpotent in the center fiber. Its characters
are therefore the \(\ell\) values
\[
 c\prod_i\rho_i\,\zeta^{\sum_i\alpha_i}.
\]
The corresponding simple modules of the full-center quotient
are exactly those with this sum; a nilpotent central element
acts by zero on every simple module.

Finally, suppose that a point above \(\chi_0\) lies in a seed
chart with unit set \(U\). Theorem~\ref{thm:intro-chart} implies
\[
 |K_U|=\ell^{|U|},\qquad
 \ell^{m-1}=\ell^{|U|-b},\qquad b\le m,
\]
where \(b\) is the number of radical supports entirely in \(U\).
Thus \(|U|\le 2m-1\), while the number of mutable variables
is \(\binom{2m}{2}-2m=2m^2-3m>2m-1\) for \(m\ge3\).
For \(m=2\) the two mutables and \(H\) force \(|U|=3\),
so \(b=2\); the second complete support requires \(F_1,B\)
to be units as well, a contradiction.
\end{proof}

\begin{example}
For \(n=4\) the full-center layers in this family are also
explicit. The crossing entries are \(x_{13},x_{14},x_{23},x_{24}\).
Here \(F_1=x_{14}\), \(B=x_{23}\); they commute with one another
and are scalar-normal in \(A\). In a block with fixed value of
\(H\), the extra central ideal is generated by
\(x_{14}^{h}x_{23}^{\ell-h}\), \(1\le h<\ell\).
It removes exactly the crossing exponent vectors with
\(h_{14}+h_{23}\ge\ell\). The \(H^\ell\) relation at this
power character has distinct roots; the central idempotent
selecting the specified root only restricts \(\sum_i\alpha_i\).
Thus, with \(\gamma=\varepsilon_1-\varepsilon_2\),
\[
 \sum_{\beta:\,\sum_i\beta_i=\sum_i\alpha_i}
       C_{\beta\alpha}(t)u^{\beta-\alpha}
 =
 \left(\sum_{h=0}^{\ell-1}t^h u^{h\gamma}\right)^2
 \sum_{h=0}^{\ell-1}(h+1)t^h u^{h\gamma}
\]
in the full-center fiber. The two free factors record
\(x_{13},x_{24}\); the factor \(h+1\) counts the splittings
\(h=h_{14}+h_{23}\).
\end{example}

\section{A comparison in subregular type \texorpdfstring{\(A_2\)}{A2}}\label{a2:section}

\subsection{Conventions for the subregular quantum-group fiber}

Let \(P=\mathbb Z\omega_1\oplus\mathbb Z\omega_2\) be the weight
lattice of \(A_2\), normalized by
\((\omega_i,\alpha_j)=\delta_{ij}\). The Cartan matrix is
\[
 \bigl((\alpha_i,\alpha_j)\bigr)_{i,j=1}^2
     =\begin{pmatrix}2&-1\\-1&2\end{pmatrix}.
\]
Throughout this section \(q\) has primitive odd order
\(\ell>3\), \(3\nmid\ell\), and \(\delta=q-q^{-1}\).
We use the weight-lattice De Concini--Kac specialization
\(U_q(\mathfrak{sl}_3)\), with generators \(E_i,F_i\) for
\(i=1,2\) and \(K_\lambda\) for \(\lambda\in P\). Its toral and
mixed relations are
\[
\begin{aligned}
 K_0&=1,& K_\lambda K_\mu&=K_{\lambda+\mu},\\
 K_\lambda E_iK_\lambda^{-1}&=q^{(\lambda,\alpha_i)}E_i,
 &K_\lambda F_iK_\lambda^{-1}&=q^{-(\lambda,\alpha_i)}F_i,\\
 [E_i,F_j]&=\delta_{ij}\frac{K_i-K_i^{-1}}{\delta},
 &K_i&=K_{\alpha_i}.
\end{aligned}
\]
For \(i\ne j\), both \(G=E\) and \(G=F\) satisfy the Serre relation
\[
 G_i^2G_j-(q+q^{-1})G_iG_jG_i+G_jG_i^2=0.
\]
There are no nilpotence relations in this specialization before
a central character is fixed.

Put \(Y=F_2F_1-qF_1F_2\), and take the matching positive PBW
root vector \(E_{12}=E_2E_1-qE_1E_2\). In the triangular PBW
basis there are three ordered root factors on each side of
the torus, corresponding to \(\alpha_1,\alpha_1+\alpha_2,\alpha_2\);
on the negative side we use \(F_1,Y,F_2\).
The \emph{Frobenius center} \(Z_{\mathrm{Fr}}\) is the central
subalgebra generated by the six \(\ell\)-th root-vector powers
and \(K_\lambda^\ell\), \(\lambda\in P\); see
\citet{DKP} for the root-of-unity central structure.
For a Frobenius character \(\chi:Z_{\mathrm{Fr}}\to\mathbb C\),
write
\[
 U_\chi=U_q(\mathfrak{sl}_3)\big/
   \bigl\langle z-\chi(z):z\in Z_{\mathrm{Fr}}\bigr\rangle.
\]
The root-vector exponents in the resulting PBW basis are
\(0,\ldots,\ell-1\), and torus exponents are taken modulo
\(\ell P\). We take
\[
  \begin{gathered}
 \chi(Y^\ell)=b\ne0,\qquad \chi(K_\lambda^\ell)=1,\\
 \chi(F_1^\ell)=\chi(F_2^\ell)
 =\chi(E_1^\ell)=\chi(E_{12}^\ell)=\chi(E_2^\ell)=0.
 \end{gathered}
\]
In type \(A_2\), the subregular unipotent class has Jordan
type \((2,1)\). This is its Frobenius representative, up to
the normalization of \(b\); see Brown--Gordon
\citep[Section 4.8]{BG} for the representative and the central
decomposition.

The \emph{Harish--Chandra center} \(Z_{\mathrm{HC}}\) is the
central subalgebra identified by the Harish--Chandra projection
with the shifted Weyl-invariant algebra on the Harish--Chandra torus; we use the two
fundamental trace elements in the normalization of
\citet[Example~2.33]{IpMan}. Let \(\mathfrak m_\eta\)
be the kernel of a compatible Harish--Chandra character.
The \emph{generalized} block consists of the \(U_\chi\)-modules
annihilated by some power of \(\mathfrak m_\eta\); the
\emph{exact} reduction consists of those annihilated by
\(\mathfrak m_\eta\) itself. Thus the exact reduction is a quotient
of the generalized block, not another name for it.
A regular Harish--Chandra orbit has trivial Weyl stabilizer in
the shifted torus coordinates. Tanisaki's Theorem~10.4
\citep{Tanisaki} concerns projectives in the generalized,
centralizer-torus-graded category, under the hypotheses of that
theorem.

Here the same signed-path technique computes projective multiplicities
in a De Concini--Kac block. The cubic exchange in this section is
not asserted to be an ordinary quantum seed mutation.

\begin{proposition}\label{a2:exchange}
There is a matrix decomposition \(U_\chi\cong M_{\ell^2}(B)\) that
takes central elements to scalar matrices over \(B\). The corner
\(B\) is generated by \(R,X,W\), where
\begin{equation}\label{a2:weights}
 R^\ell=1,\qquad RX=qXR,\qquad RW=q^{-1}WR,\qquad X^\ell=W^\ell=0.
\end{equation}
Let \(s,r\) be the compressed first and second fundamental
Harish--Chandra elements in the normalization of
\citet[Example~2.33]{IpMan}. Then
\begin{equation}\label{a2:product}
 XW=h(R),\qquad WX=h(qR),\qquad
 h(z)=\delta^{-3}(z^3+s z^{-1}-rz-z^{-3}).
\end{equation}
In particular \(\delta^3 R^3h(R)=(R^2)^3-r(R^2)^2+sR^2-1\).
These formulas hold before quotienting by the maximal ideal of the
Harish--Chandra center.
\end{proposition}

\begin{proof}
Set \(T=K_1\), \(S=K_{\omega_1-\omega_2}\), and
\(U=F_1\), \(V=Y^{-1}F_2\).
For the Serre calculation, expand \(Y\) and use the relations for
\(F_1,F_2\):
\[
\begin{aligned}
 YF_1&=F_2F_1^2-qF_1F_2F_1
      =q^{-1}F_1F_2F_1-F_1^2F_2=q^{-1}F_1Y,\\
 F_2Y&=F_2^2F_1-qF_2F_1F_2
      =q^{-1}F_2F_1F_2-F_1F_2^2=q^{-1}YF_2,\\
 VU&=Y^{-1}(Y+qF_1F_2)=1+q^2UV .
\end{aligned}
\]
Also \(U^\ell=V^\ell=0\). Since \(TY=q^{-1}YT\) and \(T^\ell=1\),
the idempotent \(e_T=\ell^{-1}\sum_{k=0}^{\ell-1}T^k\)
gives matrix units \(Y^i e_TY^{-j}\), \(0\le i,j<\ell\):
\[
  (Y^ie_TY^{-j})(Y^ke_TY^{-l})
     =\delta_{jk}Y^ie_TY^{-l},\qquad
  \sum_{i=0}^{\ell-1}Y^ie_TY^{-i}=1.
\]
Set \(x=UY^{-2}T\), \(z=q^2VY^2T^{-1}\). The exponents for
conjugation by \(Y,T\) on \(U,V\) are \((-1,1)\) and \((-2,2)\),
respectively, and \(TYT^{-1}=q^{-1}Y\). Thus
\[
\begin{aligned}
 YxY^{-1}=TxT^{-1}=x,\qquad &YzY^{-1}=TzT^{-1}=z,\\
 xz=q^2U Y^{-2} VY^2=UV,\qquad
 &zx=q^2 VY^2UY^{-2}=VU,\qquad x^\ell=z^\ell=0 .
\end{aligned}
\]
The algebra on these two elements is
\(M_\ell(\mathbb C)\): ordered monomials give the dimension bound
\(\ell^2\), and on vectors \(v_0,\ldots,v_{\ell-1}\) the operators
\[
 xv_j=v_{j+1},\qquad zv_j=[j]_{q^2}v_{j-1}
\]
where \([j]_{q^2}=(1-q^{2j})/(1-q^2)\) and out-of-range
vectors are zero. They satisfy the relations and generate
all matrix units. Indeed \(xz\) has distinct eigenvalues and the
neighboring transitions are nonzero. Put
\[
 \begin{gathered}
 H=1+(q^2-1)xz,\quad \eta=(\ell+1)/2,\quad R=SH^\eta,\\
 e_H=\ell^{-1}\sum_{k=0}^{\ell-1}H^k,\quad p=e_Te_H.
 \end{gathered}
\]
Thus \(H^\ell=1\), and \(H\) has conjugation weights \(2,-2\)
on \(x,z\). The weights of \(S\) on these elements are \(-1,1\);
it commutes with \(T,Y,H\). It follows that \(R\) commutes with the
product \(D\cong M_{\ell^2}(\mathbb C)\) of the two matrix factors.
The corner \(B=pU_\chi p\) gives the asserted matrix decomposition.
We identify the compression of \(R\) into \(B\) with \(R\) in the notation.

To describe the corner, filter by positive PBW height. The negative
Borel (including the torus) is unchanged in the associated graded;
it is generated by \(D\) and \(R\), as one recovers \(S,U,F_2\)
from the formulas above. Its PBW dimension is \(\ell^5\), so it is
\(D\otimes \mathbb C[R]/(R^\ell-1)\).
Positive root elements commute with the negative root algebra in
the associated graded, since mixed Chevalley commutators lower
positive height. If a homogeneous positive PBW monomial \(E\) has
\(T\)-weight \(q^a\), using its integer weight \(a\) set
\(D_a=Y^aH^{\eta a}\). The graded symbol of \(D_aE\) commutes
with \(D\): \(Y^a\) cancels the \(T\)-weight, and \(H^{\eta a}\)
cancels the quasi-commutation phases of \(E\) against \(x,z\).
Also \(D_aD_{a'}=D_{a+a'}\), with negative powers allowed.
The PBW basis therefore yields a filtered basis
\begin{equation}\label{a2:corner-pbw}
 pR^jD_aE p\qquad (0\le j<\ell).
\end{equation}
Here \(E\) ranges over the restricted positive PBW monomials.
In particular \(X=pD_2 E_1p\), \(W=pD_{-1}E_2p\), and \(R\)
generate the corner, since their symbols generate the associated
graded. Give \(E_1,E_2,F_1,F_2\) integer degrees \(1,-1,-1,1\),
respectively, and the torus degree zero; this is an additional
grading preserved by \(p\). The degrees in \eqref{a2:corner-pbw}
are between \(1-\ell\) and \(\ell-1\), giving \(X^\ell=W^\ell=0\).
The \(R\)-conjugation relations follow directly from \(pH=p\)
and the \(S\)-weights of \(E_i\).

It remains to compute the products. In the Chevalley conventions
above, the bold generators of Ip--Man are
\(\mathbf e_i= -\mathrm i\delta E_i\), \(\mathbf f_i= -\mathrm i\delta F_i\),
and their fractional torus element in (2.63)--(2.65) is \(S\).
The identities
\[
 F_2F_1-q^{-1}F_1F_2=YH,\qquad
 F_1F_2-q^{-1}F_2F_1=-q^{-1}Y
\]
give, in those Casimir formulas, the positive-height-two terms
\[
 \mathbf e_{21}\mathbf f_{12}
    =\delta^2(qE_1E_2-E_2E_1)YH,qquad
 \mathbf e_{12}\mathbf f_{21}
    =\delta^2(q^{-1}E_1E_2-E_2E_1)Y.
\]
In positive height two, write \(\equiv\) for equality modulo
height at most one. Since \(pH=p\), the graded PBW multiplication
of the corrected factors gives
\[
\begin{aligned}
 XW&\equiv pD_1E_1E_2p,& WX&\equiv pD_1E_2E_1p,\\
 R^{-1}s&\equiv\delta^2(qXW-WX),&
 Rr&\equiv\delta^2(q^{-1}XW-WX),\\
 R^{-1}s-Rr&\equiv\delta^3XW,&
 q^{-1}R^{-1}s-qRr&\equiv\delta^3WX .
\end{aligned}
\]
All expressions have signed degree zero. The height-one
PBW symbols have signed degrees \(1,-1\); thus only the height-zero
coefficient ring can contain the two remainders:
\[
\begin{aligned}
 \delta^3XW&=R^{-1}s-Rr+c(R),\\
 \delta^3WX&=q^{-1}R^{-1}s-qRr+d(R),\qquad
 c,d\in\mathbb C[R]/(R^\ell-1).
\end{aligned}
\]
For each \(j\), induce the one-dimensional positive Borel
module with \(E_iv=0\), \(Tv=v\), and \(Sv=q^jv\) to \(U_\chi\).
In this baby Verma module, \(v\) denotes the induced highest
vector and \(K_2v=kv\) with \(k=q^{-3j}\).
Then \(Rpv=q^jpv\) and \(pv\ne0\): the polynomial \(e_H\) in
\(UV\) has constant term one, and every positive power of \(UV\)
has positive outer \(U,V\) exponents upon negative PBW ordering.
Moreover \(Wpv=0\). Indeed \(E_2pv\) has degree \(-1\);
in the induced negative PBW basis, terms of this degree have more
\(F_1\)'s than \(F_2\)'s. In terms of the two pivot factors these
become scalar multiples of \(x^u z^wY^aT^b v\), \(u>w\),
and left compression after multiplying by \(D_{-1}\) kills them:
\(D_{-1}\) commutes with \(e_H\), and \(e_Hx=0\).
For the Casimir evaluations, \(E_1F_1v=0\) and
\(E_2F_2v=\delta^{-1}(k-k^{-1})v\). Of the four products
\(E_1E_2,E_2E_1\) acting on \(F_2F_1v,F_1F_2v\), the only one
that can be nonzero is
\(E_2E_1 F_1F_2v=\delta^{-1}(k-k^{-1})v\).
Ip--Man's formulas therefore give the scalar values
\[
 s=q^j((q^2+1)k+q^{-2}k^{-1}),\qquad
 r=q^{-j}((q^{-2}+1)k^{-1}+q^2k).
\]
As \(XWpv=0\), the preceding scalar values give
\[
 \begin{aligned}
  c(q^j)
  &= -q^{-j}s+q^jr\\
  &=-\bigl((q^2+1)k+q^{-2}k^{-1}\bigr)
     +(q^{-2}+1)k^{-1}+q^2k\\
  &=k^{-1}-k=q^{3j}-q^{-3j}.
 \end{aligned}
\]
The \(\ell\) distinct eigenvalues determine \(c(R)=R^3-R^{-3}\).
Writing \(Xf(qR)=f(R)X\), associativity gives
\[
 0=\delta^3\bigl((XW)X-X(WX)\bigr)
   =X\bigl(c(qR)-d(R)\bigr).
\]
In \eqref{a2:corner-pbw}, \(X\) times each \(R\)-weight
idempotent is nonzero. Hence \(d(R)=c(qR)\), proving
\eqref{a2:product}.
\end{proof}

Fix a regular Weyl orbit of Harish--Chandra characters above the
unipotent Frobenius value, and write \(B_{\mathrm{gen}}\) for the
corresponding direct factor of \(B\). By the center fiber-product
description and the reduced-trace splitting of Brown--Gordon
\citep[Sections 4.5 and 4.8]{BG}, \(B_{\mathrm{gen}}\) contains the
local Harish--Chandra algebra
\begin{equation}
 A=\mathbb C[\xi_1,\xi_2,\xi_3]/(e_1(\xi),e_2(\xi),e_3(\xi)).
 \label{eq:a2-coinvariant}
\end{equation}
Indeed, the \(\ell\)-power map is \'{e}tale on the torus, the source orbit
has trivial stabilizer, and logarithms linearize the \(S_3\)-action
at the target identity. Let \(\mathfrak m=(\xi_1,\xi_2,\xi_3)\).
The exact Harish--Chandra quotient is
\(B_{\mathrm{ex}}=B_{\mathrm{gen}}/\mathfrak m B_{\mathrm{gen}}\).

\begin{theorem}
\label{thm:a2-cartan}
The generalized block has three simple modules up to isomorphism. Its ordinary Cartan matrix has every entry equal to \(6\). In the exact Harish--Chandra quotient the Cartan matrix has diagonal entries \(1\) and off-diagonal entries \(2\). The latter quotient is Morita equivalent to the path algebra of a doubled triangle with immediate reversals and oriented cycles set to zero.

For a projective cover \(P_i\) and a simple \(L_j\), let \(L_{ij}(u)\) record the multiplicities of \(L_j\) in the successive radical layers of \(P_i\), starting with the top in degree zero. In the generalized block one has
\begin{equation}
\begin{aligned}
 L_{ii}(u)&=1+2u^2+2u^4+u^6,\\
 L_{ij}(u)&=u(1+2u^2+u^4)+u^2(1+u^2)\quad(i\ne j).
\end{aligned}
\label{eq:a2-loewy}
\end{equation}
In the exact quotient these polynomials become \(1\) and \(u+u^2\), respectively.
\end{theorem}

\begin{proof}
Let \(Q(t)=1+t+\cdots+t^{\ell-1}\). The positive Borel regular
module, filtered by positive PBW height, contains each torus
character with graded multiplicity \(\ell Q(t)Q(t^{-1})\).
Induction is exact by PBW freeness and gives baby Vermas with
that same signed Hilbert polynomial. The regular orbit has six
characters. The matrix factors compressed in Proposition~\ref{a2:exchange}
have signed Hilbert polynomials \(\ell^2\) and \(Q(t)Q(t^{-1})\),
respectively: the second has matrix degrees \(i-j\),
\(0\leq i,j<\ell\). Thus
\begin{equation}
 \operatorname{Hilb}(B_{\mathrm{gen}},t)
 =\frac{6\ell^2(Q(t)Q(t^{-1}))^2}{\ell^2Q(t)Q(t^{-1})}
 =6Q(t)Q(t^{-1}).
 \label{eq:a2-hilbert}
\end{equation}

Write \(e_i\) for the \(R=q^i\) idempotent, \(i\in\mathbb Z/\ell\).
Proposition~\ref{a2:exchange} gives
\[
 Xe_i=e_{i+1}X,\qquad We_i=e_{i-1}W,\qquad
 WX e_i=h(q^{i+1})e_i .
\]
The cubic in \(R^2\) has three distinct residual roots. Lift them
to \(\beta_1,\beta_2,\beta_3\) over \(A\) and denote their residues by
\(\bar\beta_1,\bar\beta_2,\bar\beta_3\). At the three bonds with
\(q^{2(i+1)}=\bar\beta_j\), the backtracking product is a unit times
\(\bar\beta_j-\beta_j\). Since
\[
 \log(\beta_j^\ell)=\ell\log(\beta_j/\bar\beta_j),
 \qquad \bar\beta_j^\ell=1,
\]
these products are unit multiples of the three logarithmic
coordinates \(\xi_j\). All other bonds are invertible. This is
a cyclic weight-orbit calculation of the kind used for generalized
Weyl algebras \citep{DGO}, with the local \(\xi_j\)-coefficients
retained rather than specialized to scalars.

For a directed monotone arc \(\gamma\) of length \(|k|<\ell\),
let \(b(\gamma)\) be the number of its nonunit bonds and
\(c_\gamma\) the product on the complementary arc. Backtracking
contracts to a coefficient, and contracting the full oriented
cycle against the reverse complementary arc gives
\(c_\gamma\gamma=0\). Table~\ref{tab:arcs} records the coefficient spaces.
Here, up to a unit, \(c_\gamma=\xi_j\xi_k=\xi_i^2\) when the complement
contains two nonunit bonds; the identities used are
\[
 \xi_1+\xi_2+\xi_3=0,\qquad
 \xi_j\xi_k=\xi_i^2\quad(\{i,j,k\}=\{1,2,3\}),\qquad
 \xi_1\xi_2\xi_3=0 .
\]
\begin{table}[htbp]
\caption{Coefficient spaces for a monotone arc.}\label{tab:arcs}
\centering
\begin{tabular}{>{$}c<{$}|>{$}c<{$}>{$}c<{$}>{$}c<{$}>{$}c<{$}}
 b(\gamma)&0&1&2&3\\ \hline
 A/(c_\gamma)&A&A/(\xi_i^2)&A/(\xi_j)&0\\
 \dim_{\mathbb C} A/(c_\gamma)&6&4&2&0 .
\end{tabular}
\end{table}
For example, with \(x=\xi_1\), \(y=\xi_2\), elimination of
\(\xi_3\) gives
\[
\begin{aligned}
 A&=\mathbb C[x,y]/(x^2+xy+y^2,\ xy(x+y)),\\
 A/(x)&=\mathbb C[y]/(y^2),\\
 A/(x^2)&=\mathbb C[x,y]/(x^2,\ xy+y^2)
          =\operatorname{span}_{\mathbb C}\{1,x,y,xy\}.
\end{aligned}
\]
This gives dimensions \(2,4\) for the two intermediate columns.
All signed degrees \(|k|\geq\ell\) vanish by the power relations.
For \(|k|<\ell\), the \(\ell\) monotone arcs of degree \(k\) contain
each of the three nonunit bonds exactly \(|k|\) times. Hence
\[
 \sum_{\deg\gamma=k}\bigl(6-2b(\gamma)\bigr)
 =6\ell-2(3|k|)=6(\ell-|k|).
\]
This is an upper bound from the displayed relations. The exchange
algebra surjects onto \(B_{\mathrm{gen}}\), and the coefficient of
\(t^k\) in \eqref{eq:a2-hilbert} attains the bound. Thus the
displayed spaces give the full signed normal form.

Contract the invertible bonds and retain one vertex in each
of the three intervals. This is a full Morita corner. At a nonunit
bond, rescale one arrow by a central unit so that its two
backtracking products equal \(\xi_i\) on the respective vertices.
The oriented cycles remain zero. Put weight one on these arrows
and weight two on \(\xi_i\). Backtracking, complementary-product,
and cycle relations are homogeneous and exhaust the relations by
the normal form. The contracted entries are therefore as follows.
See Table~\ref{tab:contracted}.
\begin{table}[htbp]
\caption{Coefficient spaces in the three-vertex Morita corner.}\label{tab:contracted}
\centering
\begin{tabular}{>{$}c<{$}|>{$}c<{$}|>{$}c<{$}|>{$}c<{$}}
\text{entry}&\text{path weight}&\text{coefficient ring}
 &\text{central Hilbert series}\\ \hline
(i,i)&0&A&1+2t+2t^2+t^3\\
(i,j),\ i\ne j&1&A/(\xi_a^2)&1+2t+t^2\\
(i,j),\ i\ne j&2&A/(\xi_b)&1+t .
\end{tabular}
\end{table}
The central Hilbert series use \(\deg\xi_i=1\), not the arrow
grading. In the latter, the degree-zero algebra is \(\mathbb C^3\);
the positive part is nilpotent and equals the Jacobson radical.
Every \(\xi_i\) is recovered from backtracking products, so the
radical is generated by arrows and its \(d\)-th power is the sum
of weights at least \(d\). This gives \eqref{eq:a2-loewy} by replacing
\(t\) by \(u^2\) and multiplying the path entries by \(u,u^2\).
The three retained idempotents are primitive. Modulo \(\mathfrak m\)
the paths still have scalar coefficients, but all backtracking
products and full oriented cycles vanish. This gives the
exact-quotient assertion.
\end{proof}

The exact quotient realizes, in this subregular type \(A_2\) case, the quantum counterpart of the no-cycle algebra anticipated in the modular work of \citet[Section~1.7]{GR}. The generalized block retains the Harish--Chandra coinvariant ring and has different projective multiplicities. Its integer grading by signed degree is distinct from the radical grading in \eqref{eq:a2-loewy}. For two selected representatives in distinct intervals, let \(d\) be their positive \(R\)-step distance and \(b\) the number of nonunit bonds crossed by that arc. The signed-degree entry is
\[
 (6-2b)t^d+2b\,t^{d-\ell},
\]
and the diagonal entry is \(6\) in degree zero.

For the highest-root subregular representative, the centralizer torus is
\[
 C=\{\operatorname{diag}(z,z^{-2},z):z\in\mathbb C^\times\},
 \qquad \alpha_1|_C=z^3,\quad \alpha_2|_C=z^{-3}.
\]
Thus signed degree \(k\) has \(C\)-weight \(3k\). At \(R=q^i\), the
test element \(K_{\alpha_1-\alpha_2}=S^3=R^3\) gives weight
\(3i\) modulo \(\ell\) in the corner, as required by
\citep[Lemma 8.3]{Tanisaki}. For chosen integer lifts \(i,j\)
the two off-diagonal occurrences therefore have shifts
\[
 3(i+d-j),\qquad 3(i+d-\ell-j)
\]
with multiplicities \(6-2b,2b\). This computes multiplicities in the
algebraic category on the representation side of
\citep[Theorem 10.4]{Tanisaki} when its hypotheses hold. It neither
identifies the lifts with the geometric basis labels nor computes
the canonical Laurent polynomials occurring before evaluation at~\(1\).
The exact quotient is not the generalized block of that theorem.

% Before submission, insert the authors' confirmed funding and
% competing-interest declarations here, as applicable to this manuscript.

\subsection*{Funding.}
This work was supported in part by the National Natural Science
Foundation of China (Grant No.~12271085) and the Natural Science
Foundation of Heilongjiang Province (Grant No.~PL2024A007).

\subsection*{Declaration of competing interest}
The authors declare no competing interests.

\subsection*{Declaration of generative AI and AI-assisted technologies in the manuscript preparation process}
The authors used ChatGPT (OpenAI, GPT-6) for exploratory mathematical
discussion, literature searches, checks of intermediate calculations
and revision of the exposition.  The authors take full
responsibility for the content of the article.

\setlength{\bibsep}{5pt}


\begin{thebibliography}{99}

\bibitem[Adaloglou et al.(2026)Adaloglou, Pasquotto, and Zanardini]{APZ}
Adaloglou, N., Pasquotto, F., Zanardini, A.:
Symplectic cohomology of quasihomogeneous \(cA_n\) singularities.
Math. Z. \textbf{312}, Article 76 (2026).
\url{https://doi.org/10.1007/s00209-026-03943-7}

\bibitem[Bavula(1993)]{Bavula}
Bavula, V.V.: Generalized Weyl algebras and their representations.
St. Petersburg Math. J. \textbf{4}, 71--92 (1993)

\bibitem[Bera and Mukherjee(2026)]{BeraMukherjee}
Bera, S., Mukherjee, S.:
On simple modules over the quantum matrix algebra at roots of unity.
J. Noncommut. Geom. (2026), published online first.
\url{https://doi.org/10.4171/JNCG/700}

\bibitem[Berenstein and Zelevinsky(2005)]{BZ}
Berenstein, A., Zelevinsky, A.: Quantum cluster algebras.
Adv. Math. \textbf{195}, 405--455 (2005).
\url{https://doi.org/10.1016/j.aim.2004.08.003}

\bibitem[Brown and Gordon(2001)]{BG}
Brown, K.A., Gordon, I.: The ramification of centres: Lie algebras
in positive characteristic and quantised enveloping algebras.
Math. Z. \textbf{238}, 733--779 (2001).
\url{https://doi.org/10.1007/s002090100274}

\bibitem[Brown and Gordon(2002)]{BG02}
Brown, K.A., Gordon, I.: The ramifications of the centres:
quantised function algebras at roots of unity.
Proc. Lond. Math. Soc. (3) \textbf{84}, 147--178 (2002).
\url{https://doi.org/10.1112/PLMS/84.1.147}

\bibitem[Brown and Macleod(2017)]{BM}
Brown, K.A., Macleod, M.J.: The Cohen Macaulay property for
noncommutative rings. Algebr. Represent. Theory \textbf{20},
1433--1465 (2017).
\url{https://doi.org/10.1007/s10468-017-9694-z}

\bibitem[Brown and Yakimov(2018)]{BY}
Brown, K.A., Yakimov, M.: Azumaya loci and discriminant ideals of
PI algebras. Adv. Math. \textbf{340}, 1219--1255 (2018).
\url{https://doi.org/10.1016/j.aim.2018.10.024}

\bibitem[De Concini et al.(1992)De Concini, Kac, and Procesi]{DKP}
De Concini, C., Kac, V.G., Procesi, C.: Quantum coadjoint action.
J. Amer. Math. Soc. \textbf{5}, 151--189 (1992).
\url{https://doi.org/10.1090/S0894-0347-1992-1124981-X}

\bibitem[Drozd et al.(1996)Drozd, Guzner, and Ovsienko]{DGO}
Drozd, Yu.A., Guzner, B.L., Ovsienko, S.A.: Weight modules over
generalized Weyl algebras. J. Algebra \textbf{184}, 491--504 (1996).
\url{https://doi.org/10.1006/jabr.1996.0270}

\bibitem[Evans and Lekili(2026)]{EL}
Evans, J., Lekili, Y.: Non-commutative crepant resolutions of
\(cA_n\) singularities via Fukaya categories.
Doc. Math. \textbf{31}, 1--26 (2026).
\url{https://doi.org/10.4171/DM/1044}

\bibitem[Goodearl and Yakimov(2017)]{GY}
Goodearl, K.R., Yakimov, M.: Quantum cluster algebra structures
on quantum nilpotent algebras. Mem. Amer. Math. Soc.
\textbf{247}, no. 1169 (2017).
\url{https://doi.org/10.1090/memo/1169}

\bibitem[Goodearl and Yakimov(2023)]{GYPoisson}
Goodearl, K.R., Yakimov, M.T.: Cluster algebra structures on Poisson
nilpotent algebras. Mem. Amer. Math. Soc. \textbf{290}, no. 1445 (2023).

\bibitem[Gordon and Rumynin(2003)]{GR}
Gordon, I., Rumynin, D.: Subregular representations of
\(\mathfrak{sl}_n\) and simple singularities of type \(A_{n-1}\).
Compos. Math. \textbf{138}, 337--360 (2003).
\url{https://doi.org/10.1023/A:1027381710548}

\bibitem[Huang et al.(2023)Huang, L\^e, and Yakimov]{HLY}
Huang, S., L\^e, T.T.Q., Yakimov, M.: Root of unity quantum
cluster algebras and Cayley--Hamilton algebras. Trans. Amer. Math.
Soc. \textbf{376}, 6767--6798 (2023).
\url{https://doi.org/10.1090/tran/8904}

\bibitem[Ip and Man(2024)]{IpMan}
Ip, I.C.-H., Man, R.: Positive representations with zero Casimirs.
Comm. Math. Phys. \textbf{405}, article 91 (2024).
\url{https://doi.org/10.1007/s00220-024-04955-2}

\bibitem[Iyama and Wemyss(2014)]{IW}
Iyama, O., Wemyss, M.: Singular derived categories of
\(\mathbb Q\)-factorial terminalizations and maximal modification
algebras. Adv. Math. \textbf{261}, 85--121 (2014).
\url{https://doi.org/10.1016/j.aim.2014.04.001}

\bibitem[Jing and Ray(2012)]{JR}
Jing, N., Ray, R.: Zonal polynomials and quantum antisymmetric
matrices. Bull. Inst. Math. Acad. Sin. (N.S.) \textbf{7}, 1--31 (2012)

\bibitem[Jing and Zhang(2024)]{JZ}
Jing, N., Zhang, J.: Minor identities for Sklyanin determinants.
Adv. Math. \textbf{442}, 109569 (2024).
\url{https://doi.org/10.1016/j.aim.2024.109569}

\bibitem[Muller et al.(2024)Muller, Nguyen, Trampel, and Yakimov]{MNTY}
Muller, G., Nguyen, B., Trampel, K., Yakimov, M.: Poisson geometry
and Azumaya loci of cluster algebras. Adv. Math. \textbf{453},
109822 (2024). \url{https://doi.org/10.1016/j.aim.2024.109822}

\bibitem[Nguyen et al.(2025)Nguyen, Trampel, and Yakimov]{NTY}
Nguyen, B., Trampel, K., Yakimov, M.: Root of unity quantum cluster
algebras and discriminants. J. Lond. Math. Soc. \textbf{111},
e70060 (2025). \url{https://doi.org/10.1112/jlms.70060}

\bibitem[Strickland(1996)]{Strickland}
Strickland, E.: Classical invariant theory for the quantum
symplectic group. Adv. Math. \textbf{123}, 78--90 (1996).
\url{https://doi.org/10.1006/aima.1996.0067}

\bibitem[Tanisaki(2027)]{Tanisaki}
Tanisaki, T.: Quantized flag manifolds and non-restricted modules
over quantum groups at roots of unity. J. Amer. Math. Soc.
\textbf{40}, 129--193 (2027).
\url{https://doi.org/10.1090/jams/1078}

\bibitem[Zhang and Tang(2026)]{ZT}
Zhang, Y., Tang, X.: Cluster structure on the quantum coordinate
ring of skew-symmetric matrices in general case.
J. Algebra \textbf{690}, 677--700 (2026).
\url{https://doi.org/10.1016/j.jalgebra.2025.11.016}

\end{thebibliography}
\end{document}